\documentclass[a4paper,11pt]{article}
\usepackage[top=3cm,bottom=3cm,left=2.5cm,right=2.5cm,a4paper]{geometry}
\usepackage{amsfonts,amsmath,amssymb,amsthm,upgreek}
\usepackage{mathbbol}
\usepackage{float,fancyhdr,siunitx}
\usepackage[dvipsnames]{xcolor} 
\usepackage{setspace}
\usepackage{comment}
\usepackage{lineno}
\usepackage[normalem]{ulem} 
\usepackage[colorlinks]{hyperref}
\hypersetup{
    colorlinks=true, 
    linkcolor={blue!50!black},
    citecolor={green!60!black},
    urlcolor= {red!80!black}
}
\usepackage{stmaryrd} 
\usepackage{enumitem,mdframed}
\usepackage{empheq,longtable}
\usepackage[bottom]{footmisc} 
\usepackage{tikz,pgfplots}
\usepgfplotslibrary{groupplots}
\usetikzlibrary{positioning,arrows,fit,patterns}
\pgfplotsset{compat=newest}
\usepackage{subfig}
\DeclareCaptionLabelFormat{opening}{\textbf{#2)}} 
\makeatletter
\renewcommand*\env@matrix[1][\arraystretch]{%
  \edef\arraystretch{#1}%
  \hskip -\arraycolsep
  \let\@ifnextchar\new@ifnextchar
  \array{*\c@MaxMatrixCols c}}
\makeatother

\usepackage[]{todonotes}

\newcommand*\samethanks[1][\value{footnote}]{\footnotemark[#1]}

\newcommand{\ii}          {\mathrm{i}}
\newcommand{\dd}          {\mathrm{d}}
\newcommand{\bx}          {\mathbf{x}}

\newcommand{\bv}          {\boldsymbol{v}}
\newcommand{\bw}          {\boldsymbol{w}}

\newcommand{\bxi}         {\boldsymbol{\xi}}
\newcommand{\bu}          {\boldsymbol{u}}
\newcommand{\bpsi}        {\boldsymbol{\psi}}

\newcommand{\n}           {\boldsymbol{n}} 

\newcommand{\rmax}        {{r_{\max}}}
\newcommand{\Lliouville}  {\mathcal{L}_{\mathrm{Liouville}}}
\newcommand{\Lliouvillec} {\mathcal{L}_{\mathrm{Liouville-c}}}
\newcommand{\Loriginaldiv}{\mathcal{L}_{\mathrm{original}}}

\newcommand{\mesh}{\mathcal{T}_h}
\newcommand{\ncell}{n_{\mathrm{cell}}}
\newcommand{\nface}{n_{\mathrm{face}}}
\newcommand{\face}{\mathfrak{f}}
\newcommand{\setface}{\Sigma}
\newcommand{\domain}{\mathbb{B}} 
\newcommand{\ddomain}{\partial \mathbb{B}}
\newcommand{\hdgA}{\mathbb{A}}
\newcommand{\hdgC}{\mathbb{C}}
\newcommand{\hdgK}{\mathbb{K}}
\newcommand{\hdgR}{\mathbb{R}}

\newcommand{\hdgF}{\mathbb{F}}
\newcommand{\hdgS}{\mathbb{S}}
\newcommand{\hdgB}{\mathbb{B}}
\newcommand{\hdgL}{\mathbb{L}}
\newcommand{\ndof}{n_{\mathrm{dof}}}
\newcommand{\ndofface}{n_{\mathrm{dof}}^\mathfrak{f}}

\newcommand{\zsrc}{z_{\mathrm{src}}}

\newlength{\cbarwidth}
\newlength{\cbarheight}

\newcommand{\adastra}{\texttt{Adastra}}
\newcommand{\mps}    {\texttt{Swan}}

\newcommand{\phifull}{\nabla\phi_0} 

\newcommand{\pressE}{\delta_p^{\mathrm{E}}}
\newcommand{\err}{\mathfrak{e}}
\newcommand{\Amat}{\boldsymbol{A}}

\newcommand{\meshA}{$\texttt{mesh06\textsc{M}}$}
\newcommand{\meshB}{$\texttt{mesh08\textsc{M}}$}
\newcommand{\meshC}{$\texttt{mesh11\textsc{M}}$}

\newcommand{\blr}[1]{\texttt{blr-}\num{e-{#1}}}
\newcommand{\fr}{FR}
\newcommand{\epsilonblr}{\varepsilon_{\texttt{BLR}}}
\newcommand{\pctNop}{$\texttt{N}_{\texttt{op}}^{\si{\percent}}$}
\newcommand{\pctNen}{$\texttt{N}_{\texttt{entries}}^{\si{\percent}}$}
\newcommand{\pctNMP}{$\texttt{N}_{\texttt{entries}}^{\si{\percent}\texttt{MP}}$}

\newcommand{\mumps}{\texttt{MUMPS}}
\newcommand{\haven}{\texttt{Hawen}}
\DeclareSymbolFont{yhlargesymbols}{OMX}{yhex}{m}{n} 
\DeclareMathAccent{\yhwidehat}{\mathord}{yhlargesymbols}{"62}

\newtheorem{remark}     {Remark}

\usepackage[capitalize,nameinlink]{cleveref}
\crefname{section}   {Section}   {Sections}
\crefname{subsection}{Subsection}{Subsections}
\Crefname{section}   {Section}   {Sections}
\Crefname{subsection}{Subsection}{Subsections}
\Crefname{figure}    {Figure}    {Figures}

\crefformat{equation}{\textup{#2(#1)#3}}
\crefrangeformat{equation}{\textup{#3(#1)#4--#5(#2)#6}}
\crefmultiformat{equation}{\textup{#2(#1)#3}}{ and \textup{#2(#1)#3}}
{, \textup{#2(#1)#3}}{, and \textup{#2(#1)#3}}
\crefrangemultiformat{equation}{\textup{#3(#1)#4--#5(#2)#6}}%
{ and \textup{#3(#1)#4--#5(#2)#6}}{, \textup{#3(#1)#4--#5(#2)#6}}{, and \textup{#3(#1)#4--#5(#2)#6}}
\Crefformat{equation}{#2Equation~\textup{(#1)}#3}
\Crefrangeformat{equation}{Equations~\textup{#3(#1)#4--#5(#2)#6}}
\Crefmultiformat{equation}{Equations~\textup{#2(#1)#3}}{ and \textup{#2(#1)#3}}
{, \textup{#2(#1)#3}}{, and \textup{#2(#1)#3}}
\Crefrangemultiformat{equation}{Equations~\textup{#3(#1)#4--#5(#2)#6}}%
{ and \textup{#3(#1)#4--#5(#2)#6}}{, \textup{#3(#1)#4--#5(#2)#6}}{, and \textup{#3(#1)#4--#5(#2)#6}}

\crefdefaultlabelformat{#2\textup{#1}#3}

\crefname{proposition}{Proposition}{Propositions}
\Crefname{proposition}{Proposition}{Propositions}
\crefname{definition} {Definition} {Definitions}
\Crefname{definition} {Definition} {Definitions}
\crefname{theorem}    {Theorem}    {Theorems}
\Crefname{theorem}    {Theorem}    {Theorems}
\crefname{remark}     {Remark}     {Remarks}
\Crefname{remark}     {Remark}     {Remarks}
\crefname{assumption} {Assumption} {Assumptions}
\Crefname{assumption} {Assumption} {Assumptions}

\title{3D Modeling of Solar Oscillations
       with Hybridizable Discontinuous Galerkin Method}

\author{
Florian Faucher\thanks{Project-Team Makutu, Inria, University of Pau and Pays de l'Adour, TotalEnergies, CNRS UMR 5142, France, (\href{mailto:florian.faucher@inria.fr}{\texttt{florian.faucher@inria.fr}}).} 
\and Ha Pham\samethanks[1]
\and Damien Fournier\thanks{Max-Planck-Institut f\"ur 
                      Sonnensystemforschung, Justus-von-Liebig-Weg 3, 
                      G\"ottingen, Germany.}
\and Patrick Amestoy\thanks{Mumps Technologies, France.}
\and H\'el\`ene Barucq\samethanks[1]
\and Jean-Yves  L'Excellent\samethanks[3]
\and Th\'eo Mary\thanks{Sorbonne University, CNRS, LIP6, France.}
\and Laurent Gizon\samethanks[2]
                 \thanks{Institut f\"ur Astrophysik und Geophysik,
                         Georg-August-Universit\"at G\"ottingen, Germany.}
}

\date{\today}

\begin{document}
\maketitle 

\begin{abstract}

With increasing quantity and quality of solar observations, 
it becomes essential to account for three-dimensional 
heterogeneities in wave modeling for seismic data interpretation.
In this context, we present a 3D solver of the time-harmonic 
adiabatic stellar oscillation equations without background flows 
on a domain consisting of the Sun and its photosphere. 
The background medium consists of 3D heterogeneities on top of a 
radial strongly-stratified standard solar model. 
The oscillation equations are solved with the 
Hybridizable Discontinuous Galerkin (HDG) 
method, considering a first-order formulation
in terms of the vector displacement and the pressure perturbation.
This method combines the high-order accuracy and the parallelism
of DG methods while yielding smaller linear systems. 
These are solved with a direct solver, 
with block low-rank compression and mixed-precision arithmetic
to reduce memory footprint.
The trade-off between compression and solution accuracy is investigated,
and our 3D solver is validated by comparing with resolution 
under axial symmetry for solar backgrounds.
The capacity of the solver is illustrated with wave speed heterogeneities 
characteristic of two physical phenomena: active regions and convection. 
We show the importance of global 3D gravito-acoustic wave simulations, 
in particular when the amplitudes of the perturbations are strong and 
their effect on the wavefield cannot be estimated by linear approximations.

\end{abstract}
\tableofcontents

\numberwithin{equation}{section}
\section{Introduction}

The Sun is a highly heterogeneous 3D rotating plasma \cite{Stix2002},
in which different types of waves are present and have been employed
to infer the properties of solar interior, 
starting from the discovery of the 5-mins solar oscillations \cite{Leighton1962} 
to the recent observations of inertial waves \cite{Loeptien2018,Gizon2021}.
Heterogeneities can be due for example to turbulent convection, responsible for 
the excitation of waves, or the magnetic activities in the Sun.
With increasing resolution in solar observations, 
accounting for three-dimensional heterogeneities 
is essential to fully interpret helioseismic data. 
This requires a full three-dimensional solver of the 
stellar oscillation equations, capable of handling
complex heterogeneities on top of the standard 
strongly-stratified solar background model (\cite{christensen1996current}), 
while mitigating the computational cost. 
Such perturbations also imply moving beyond the current assumptions of 
spherical symmetry \cite{christensen2014lecture,Pham2020Siam,Pham2024assembling}
and axial symmetry \cite{Gizon2017,Pham2025axisymreport},
as well as beyond the Born approximation \cite{Gizon2010,Svanda2011} 
which becomes inaccurate in case of strong perturbations, cf. \cite{Jackiewicz2007}.

We recall that symmetry assumptions are employed in literature for dimension 
reduction and thus lower computational cost. 
For spherical symmetry, the equations of stellar oscillations, via vector
spherical harmonics, can be decomposed to a set of decoupled 1D systems, which
have been studied in terms of eigenmodes and Green's functions, 
\cite{unno1979nonradial,christensen2014lecture,Pham2020Siam,Pham2021Galbrun,Pham2024assembling}.
Eigenvalue solvers have been developed for rapidly rotating stars \cite{Reese2006} 
and for the Sun to investigate inertial modes which are restored by the Coriolis force, 
\cite{Bekki2022}. These require axisymmetric formulations to account for the effects 
of differential rotation. Even in 2D, computations can already be expensive
due to the strong solar stratification in the near-surface layer (\cite{christensen1996current}); 
for this reason, computations in \cite{Bekki2022} are restricted to a shell 
up to $\num{0.985}R_\odot$, with $R_\odot$ the solar 
radius. Since this near-surface layer contains important
physics for seismology, in \cite{Pham2025axisymreport}, 
solar Green's kernels for axisymmetric backgrounds are 
computed for domains extending to $\num{1.001}R_\odot$ in the photosphere.

We consider the equations of stellar oscillations \cite{LyndenBell1967,unno1979nonradial,Gough1990,christensen2014lecture}
without flow with gravity prescribed in the Cowling's approximation. 
The equations are solved in a first-order formulation
in terms of the displacement and the Eulerian perturbation of 
pressure as proposed in \cite{Pham2025axisymreport};
we refer to \cite{Pham2025axisymreport} for a 
discussion of the different variants of these equations.
The challenge of modeling realistic solar (or more generally, stellar) 
backgrounds is twofold. Firstly, the model parameters vary drastically 
near the surface, for instance the solar density decreases exponentially 
in the atmosphere, \cite{christensen1996current,AtmoI2020}, see \cref{section:implementation}. 
This behaviour leads to short wavelengths and requires refined spatial 
discretization. Secondly, the inclusion of gravity effects introduces
a non-zero buoyancy frequency squared $N(\bx)^2$, \cite[Eq. (3.27)]{Pham2025axisymreport},
whose sign crucially impacts the problem. 

When $N^2$ is negative, as it is for the Sun below surface, 
convective instabilities arise \cite{christensen2014lecture}, 
rendering time-domain approaches challenging, \cite{Schunker2011,Papini2014}.
This motivates our use of the time-harmonic formulation.
With $\omega$ the frequency, the sign of $(\omega^2-N^2)$ further 
plays a key role in determining the nature of the partial differential 
equations (PDE) at zero attenuation. 
Although attenuation is always included in simulations, this change still 
requires careful treatment, cf. \cite{Pham2025axisymreport}.
In the solar case, all frequencies below the acoustic cut-off 
($\sim \num{5.3}\si{\milli\Hz}$) exhibit a coexistence of 
elliptic and hyperbolic regions, \cite{Pham2025axisymreport}. 
In applications, the frequency range of interest lies
between \num{1.5} and \num{5}\si{\milli\Hz}, 
(see for example a solar power spectrum in \cite{Gizon2010}),
thus accurate global wave simulations must be able to 
handle both the change in the nature of the PDE and the 
regions where $N^2 < 0$, \cite{Pham2025axisymreport}. 
Note that neglecting gravity effects, \cite{Gizon2017,Pham2020Siam}, 
corresponds to setting $N^2 = 0$, thereby avoiding complications. 
While this approximation can reproduce acoustic modes with a 
reasonable degree of accuracy, incorporating a non-zero $N^2$ 
yields a more physically faithful representation of stellar 
oscillations, as demonstrated in the comparisons presented 
in \cite{Pham2024assembling,Pham2025axisymreport}.


For the discretization of the oscillation equations, 
we employ the Hybridizable Discontinuous Galerkin (HDG) method, which has 
been popularized by the pioneering works of Cockburn and collaborators, e.g., 
\cite{Arnold2002,Cockburn2009,nguyen2009implicit,cockburn2010hybridizable,
      cockburn2016static,cockburn2023hybridizable}. 
The HDG discretization with $hp$-adaptivity 
is implemented in the open-source software \haven~\cite{Hawen2021}, 
and for the solution of the resulting linear systems, 
we employ the direct solver \mumps, \cite{Amestoy2001,ablm:17b}.
The HDG method works with a first-order system of equations, introducing 
the numerical trace as an additional unknown variable alongside the 
volume unknowns. 
This approach allows for the global system to be formulated in terms 
of the numerical trace only, leading to a reduction in the number 
of unknowns by eliminating the interior degrees of freedom, a feature 
which is critical for large-scale applications, e.g.,
\cite{Pham2024stabilization,Faucher2020adjoint}. 
The HDG method preserves the flexibility and advantages 
of the Discontinuous Galerkin (DG) family, including 
$hp$-adaptivity and parallelizability, while 
reducing the size of the global linear system by considering the 
numerical trace only. In addition, it provides both unknowns
of the first-order system (in our case the displacement and 
pressure perturbation) without requiring any post-processing.
The balance of efficiency and flexibility makes HDG particularly 
attractive for complex problems, we refer to \cite{Pham2024stabilization,cockburn2016static,cockburn2023hybridizable}
for additional discussion and references, and to \cite{Kirby2012,Kirby2016}
for comparisons with the Continuous Galerkin method.
Concerning the singularity propagation phenomena in the near-surface 
layer related to $N^2> 0$, we employ the stabilization proposed 
in \cite{Pham2025axisymreport} that has shown to retain the robustness 
of the method with small attenuation.

The choice of a direct solver is motivated by two factors. 
Firstly, it allows efficient computation for multiple 
right-hand sides (i.e., multiple sources of wave excitation), 
which is critical for many helioseismic applications that work
with the Green's functions or use the cross-correlation, 
cf.~\cite{Pham2020Siam,Pham2024assembling}. 
In particular the cross-correlation can be performed between
any points of observed solar disks which, with 
full resolution HMI (Helioseismic and Magnetic Imager)
\cite{Scherrer2012} 
dopplergrams means \num{4096}$^2$ observation points, 
i.e., \num{4096}$^2$ possible right-hand sides in the simulation 
to fully reproduce the observables.
Secondly, iterative solvers would be an alternative but depend 
heavily on the robustness of their preconditioners, which are 
challenging to design for solar backgrounds with strong 
variations \cite{Preuss2021}.

The main drawback of a direct solver (compared to an iterative one) 
is the complexity (number of operations and memory requirement) of the 
matrix factorization. 
In many applications (such as those coming from the discretization of elliptic PDE), 
the linear system has been shown to have a low-rank property: conveniently 
defined off-diagonal blocks can be approximated by low-rank products \cite{bebe:08} 
at a target accuracy provided by the user. 
This enables the user to control the requested accuracy and to reduce the 
complexity of the linear solver.
In this work we employ block low-rank (BLR) compression \cite{aabblw:15,ablm:17b}
and combine it with mixed-precision techniques \cite{Amestoy2023mixed}. 
Our approach to mixed precision is quite generic:
given the accuracy at which the BLR compression is performed, the mixed-precision algorithm automatically 
adapts the representation of the data (in our case the matrix of factors) 
to a continuum of precisions (up to seven precisions) to reduce the memory footprint
without further altering the accuracy of the data. 
We also investigate in the article the trade-off between 
compression and solution accuracy.  
Furthermore, the special structure of the HDG global 
matrix which involves only the degrees of freedom on the faces 
of elements to design a compact representation of the structure of the 
matrix is used to efficiently perform all the symbolic steps of the
linear solver.
We also refer to \cite{oaabbcd:23,oabb_EAGE:24} that use 
\mumps~in a different numerical context 
(acoustic wave equation with finite difference discretization and lower precision requirement)
for 3D frequency-domain full-waveform inversion. 

Our work presents the first results of fully three-dimensional 
stellar gravito-acoustic wave propagation simulations in the Sun, performed with 
a careful managing of computational cost. 
We model the Sun with both the interior and the atmosphere up to
the photosphere, for frequencies up to $2\si{\milli\Hz}$  using 
realistic background profile with gravity (model S, \cite{christensen1996current}), 
and superimposed 3D heterogeneities represented as perturbations 
in wave speed. 
We first validate our solver by comparing the 3D solution to the one obtained 
in 2.5D assuming a spherically symmetric backgrounds \cite{Pham2025axisymreport}.
Two 3D test cases are considered:
first, we use perturbations representing a thermal sunspot, 
that we create by converting the magnetic map observed at 
the surface of the Sun from GONG\footnote{\url{https://nso.edu/telescopes/nisp/gong/}} 
(Global Oscillations Network Group) observations.
It contains the details of the spot such as the umbra and 
penumbra regions, and is converted into a wave speed perturbation. 
Second, we use the 3D sound speed from a snapshot of a nonlinear convection simulation, \cite{Noraz2025}. 
This 3D convection simulation solves the nonlinear equations of momentum, continuity, and energy in a 
spherical shell from the base of the convection zone up to $0.99 R_\odot$. Due to the absence of the surface, 
the traditional pressure modes are not excited  by granulation. We can however use our 3D solver to study 
the impact of the medium obtained from the nonlinear simulation on the acoustic waves.

The paper is organized as follows.
In \cref{section:wave-eq}, we present the wave equations for 
stellar oscillations, 
and \cref{section:hdg-discretization} details the Hybridizable 
Discontinuous Galerkin (HDG) discretization method.
The features implemented to assemble and solve the resulting 
linear systems are discussed in \cref{section:implementation,section:mumps}, 
where we also investigate the trade-offs between solution accuracy 
and computational efficiency.
The validity of the 3D solver for a spherically symmetric solar background is 
established in \cref{section:validation-with-axisym} via comparisons with 
the axisymmetric solver developed in \cite{Pham2025axisymreport}.
Finally, \cref{section:numerics-3d} presents simulations of solar 
wave propagation in the Sun with 3D heterogeneities modeled as 
wave speed perturbations.

\section{Time-harmonic stellar oscillation equations}
\label{section:wave-eq}

The linear adiabatic stellar equations describe small oscillations
on top of a time-independent adiabatic background in hydrostatic equilibrium.
We specify below the formulations employed in this work, 
which were proposed in our previous work \cite{Pham2025axisymreport}.
We also refer therein for detailed discussion of derivation and 
literature of other variants existing in literature.

\subsection{Original system of equation}
\label{subsection:original-equation}
%
%
We denote by $\domain$ the sphere centered at the origin of the Sun 
with radius $\rmax$, and by $\partial \mathbb{B}$ its boundary.
Following \cite{Pham2025axisymreport}, we consider the 
first-order system of equations in terms of the Eulerian pressure
$\pressE$ and the vector-field displacement $\bxi$, which are solutions to,
\begin{subequations} \label{main-equation:original-div_v0}
\begin{empheq}[left={\empheqlbrace}]{align}
   &    \rho_0 \bigg[ - \omega^2    \,-\,2 \ii\omega\gamma_{\mathrm{att}}
      \,+\, \phifull \otimes
      \Big(- \dfrac{\nabla\rho_0}{\rho_0} \,+\, \dfrac{\nabla p_0}{\rho_0 c_0^2}
      \Big) 
      \bigg] \bxi
  \,+\, \dfrac{\phifull}{c_0^2}  \, \pressE  
  \,+\, \nabla \pressE 
  \,=\, 0 \,, \quad \text{in $\domain$}\,,\\ 
  & \dfrac{\pressE}{\rho_0 c_0^2} 
    \,+\, \bxi \cdot \dfrac{\nabla p_0}{\rho_0 c_0^2}
\,+\, \nabla\cdot\bxi \,=\, f \,,\hspace*{20.8em} \text{in $\domain$}. 
\end{empheq} \end{subequations}
Here, $f\in L^2(\domain)$ is the interior source term, $\omega\in \mathbb{R}_{> 0}$ is the frequency, 
and attenuation is incorporated with $\gamma_{\mathrm{att}}\in \mathbb{R}_{\geq 0}$. 
The background 
is characterized by the density $\rho_0$, the pressure $p_0$, the 
gravitational potential $\phi_0$, adiabatic sound speed $c_0$, and 
adiabatic index $\Gamma_1$, which satisfy the following constraints:
\begin{equation}\label{main:constraints}
\begin{aligned}
     & \Delta \phi_0 = 4\pi G\rho_0 \,, 
     \qquad \text{ with $\phi_0 \in L^2(\mathbb{R}^3)$ and $G$ the gravitational constant;}  \\
     & c_0^2 \, \rho_0 \,=\, \Gamma_1 \, p_0\, \quad \text{(Adiabacity)}\,, 
     \hspace*{3em} \nabla p_0 \,=\, - \rho_0 \, \nabla\Phi_0 \quad \text{(Hydrostatic equilibrium)}\,.
\end{aligned}\end{equation}

In the following, we introduce the scale height vectors associated to
$\rho_0$ and $c_0$
\begin{equation}
\boldsymbol{\alpha}_\rho := - \dfrac{\nabla\rho_0}{\rho_0}, \quad \boldsymbol{\alpha}_c := - \dfrac{\nabla c_0}{c_0}. 
\end{equation}
Employing the constraints \cref{main:constraints} 
to replace $\nabla p_0$ in \cref{main-equation:original-div_v0}, 
the resulting equations only contain as 
background parameters $\rho_0$, $\phi_0$ and $c_0$: 
\begin{subequations} \label{main-equation:original-div}
\begin{empheq}[left={\empheqlbrace}]{align}
   &    \rho_0 \bigg[ - \omega^2    \,-\,2 \ii\omega\gamma_{\mathrm{att}}
      \,+\, \phifull \otimes
      \Big(\boldsymbol{\alpha}_\rho \,-\, \dfrac{\nabla \phi_0}{c_0^2}
      \Big) 
      \bigg] \bxi
  \,+\, \dfrac{\phifull}{c_0^2}  \, \pressE  
  \,+\, \nabla \pressE 
  \,=\, 0 \,, \qquad \text{in $\domain$}\,,\\ 
  & \dfrac{\pressE}{\rho_0 c_0^2} 
    \,-\, \bxi \cdot \dfrac{\nabla \phi_0}{c_0^2}
\,+\, \nabla\cdot\bxi \,=\, f \,,\hspace*{20em} \text{in $\domain$}. 
\end{empheq} \end{subequations}
The above system is closed by imposing vanishing Lagrangian pressure perturbation
on the boundary $\ddomain$. It corresponds to
\begin{equation} \label{eq:boundary-condition:original}
  \pressE \,-\, \rho_0 \bxi \cdot \nabla \phi_0 
  \,=\, 0\,, \qquad \text{vacuum boundary condition on $\ddomain$.}
\end{equation}


\subsection{Liouville variants}
\label{subsection:liouville-equations}

To handle the exponential decay of $\rho$ in the near surface
layer and low atmosphere which can lead to numerical difficulties \cite{Pham2020Siam, Pham2024assembling},
  we next eliminate the presence 
of $\rho_0$ in \cref{main-equation:original-div}
by employing the two changes of variables, proposed in \cite[Section~3.2]{Pham2025axisymreport},
\begin{subequations}\begin{align}
&\text{Liouville unknowns:}  \hspace*{3.7em} \bu_L \,=\, \sqrt{\rho_0} \, \bxi \,,\hspace*{2.2em}
                             w_L \,=\, \dfrac{1}{\sqrt{\rho_0}} \pressE \,, \\
&\text{Liouville-c unknowns:} \hspace*{3em}
       \bu_c \,=\, c_0 \sqrt{\rho_0} \, \bxi \,, \hspace*{1.5em} 
  w_c \,=\, \dfrac{1}{c_0 \sqrt{\rho_0}} \pressE \,.
\end{align} \end{subequations}
The boundary value problem \cref{main-equation:original-div}-\cref{eq:boundary-condition:original}
is equivalent to the \emph{Liouville} formulation, written in terms of $(\bu_L, w_L)$,
\begin{subequations} \label{main-equation:Liouville}
\begin{empheq}[left={\empheqlbrace}]{align}
   &  \bigg[ - \omega^2    \,-\,2 \ii\omega\gamma_{\mathrm{att}}
      \,+\, \phifull \otimes
      \Big(\boldsymbol{\alpha}_\rho \,-\, \dfrac{\nabla \phi_0}{c_0^2}
      \Big) 
      \bigg] \bu_{\mathrm{L}}
  \,+\, \left( \dfrac{\phifull}{c_0^2} \,-\,\dfrac{\boldsymbol{\alpha}_\rho}{2} \right)
        \, w_{\mathrm{L}}
  \,+\, \nabla w_{\mathrm{L}}
  \,=\, 0 \,, \\
  & \dfrac{w_{\mathrm{L}}}{c_0^2} 
    \,+\, \bu_{\mathrm{L}} \cdot 
    \left( -\dfrac{\nabla \phi_0}{c_0^2} \,+\,\dfrac{\boldsymbol{\alpha}_\rho}{2} \right)
\,+\, \nabla\cdot\bu_{\mathrm{L}} \,=\, \sqrt{\rho_0} \, f \,,\hspace*{5em} \text{in $\domain$;}\\
& w_{\mathrm{L}} \,-\, \bu_{\mathrm{L}} \cdot \nabla \phi_0 \,=\, 0\,, \qquad \text{on $\ddomain$,} \label{eq:boundary-condition:liouville}
\end{empheq} \end{subequations}
and  the \emph{Liouville-c} formulation, 
with unknowns $(\bu_c, w_c)$,
\begin{subequations} \label{main-equation:Liouville-c}
\begin{empheq}[left={\empheqlbrace}]{align}
   &  \bigg[ \dfrac{-\omega^2  - 2 \ii\omega\gamma_{\mathrm{att}}}{c_0^2}
      +\dfrac{\phifull}{c_0^2} \otimes
      \Big(\boldsymbol{\alpha}_\rho - \dfrac{\nabla \phi_0}{c_0^2}
      \Big) 
      \bigg] \bu_{\mathrm{c}}
  + \left( \dfrac{\phifull}{c_0^2} -\dfrac{\boldsymbol{\alpha}_\rho}{2}  
              -\boldsymbol{\alpha}_{c} \right)
        w_{\mathrm{c}}
  + \nabla w_{\mathrm{c}}
  = 0 , \\
  & w_{\mathrm{c}}
    \,+\, \bu_{\mathrm{c}} \cdot 
    \big( -\dfrac{\nabla \phi_0}{c_0^2} \,+\,\dfrac{\boldsymbol{\alpha}_\rho}{2} 
                                           \,+\,\boldsymbol{\alpha}_{c}\big)
\,+\, \nabla\cdot\bu_{\mathrm{c}} \,=\, c_0\,\sqrt{\rho_0} \, f ,
\hspace*{5em} \text{in $\domain$;}\\
 & w_{\mathrm{c}} \,-\, \bu_{\mathrm{c}} \cdot \dfrac{\nabla \phi_0}{c_0^2} 
                 \,=\, 0\,, \qquad \text{on $\ddomain$.}
\end{empheq} \end{subequations}

\paragraph{Generic notation}
The systems of equations \cref{main-equation:original-div,main-equation:Liouville,main-equation:Liouville-c}
can be unified under the following form: find $(\bu_\bullet,w_\bullet)$ that solve
\begin{subequations}\label{eq:main-generic}
\begin{empheq}[left={\empheqlbrace}]{align}
&  \boldsymbol{A}_\bullet \, \bu_\bullet \,+\, \boldsymbol{\beta}_{\bullet} \, w_\bullet \,+\, \nabla w_\bullet \,=\, 0, 
\hspace*{4.7em} \text{in $\domain$,} \\
&  \nabla\cdot    \, \bu_\bullet \,-\, \boldsymbol{\beta}_{\bullet} \,\cdot\,\bu_\bullet \,+\,\varrho_\bullet\, w_\bullet \,=\, f_\bullet, \hspace*{2.9em} \text{in $\domain$,} \\
& w_\bullet \,+\, \mathbf{Z}_{\mathrm{bc}\bullet}\,\cdot\, \bu_\bullet \,=\, 0 \,,   \hspace*{7.9em}
    \text{on $\ddomain$,}
\end{empheq}\end{subequations}
with corresponding coefficients functions $\Amat_\bullet$, $\boldsymbol{\beta}_\bullet$ and $\varrho_\bullet$. 
The index $_\bullet$ serves to identify the variants, and we respectively use 
indices $\mathrm{o}$, $\mathrm{L}$ and $\mathrm{c}$ for 
\cref{main-equation:original-div,main-equation:Liouville,main-equation:Liouville-c}.
These systems are referred to as,
\begin{equation} \label{eq:equation-system-notation}
  \Loriginaldiv \begin{pmatrix}
                   \bu_{\mathrm{o}} \\
                   w_{\mathrm{o}}
                \end{pmatrix} \,=\, \begin{pmatrix}
                    \mathbf{0} \\
                    f_{\mathrm{o}} 
                \end{pmatrix} \,, \quad 
  \Lliouville \begin{pmatrix}
                   \bu_{\mathrm{L}} \\
                   w_{\mathrm{L}}
                \end{pmatrix} \,=\, \begin{pmatrix}
                    \mathbf{0} \\
                    f_{\mathrm{L}} 
                \end{pmatrix} \,, \quad
  \Lliouvillec \begin{pmatrix}
                   \bu_{\mathrm{c}} \\
                   w_{\mathrm{c}}
                \end{pmatrix} \,=\, \begin{pmatrix}
                    \mathbf{0} \\
                    f_{\mathrm{c}} 
                \end{pmatrix} \, .
\end{equation}
    
\section{Hybridizable Discontinuous Galerkin discretization} 
\label{section:hdg-discretization}

In this section we use the Hybridizable Discontinuous Galerkin 
(HDG) method to discretize the working equation \cref{eq:main-generic}.
We follow the notation and structure of 
\cite{Pham2024stabilization,Pham2025axisymreport,Faucher2020adjoint}.  
In particular, in \cite[Section~5]{Pham2025axisymreport} we work with 
a similar equation but assuming axially symmetric backgrounds 
with 2.5D HDG implementation, while now we consider the 3D problem.

\subsection{Preliminaries and notations}

\paragraph{Domain discretization}

The domain $\domain$ is discretized with non-overlapping cells (also referred to as elements)
that form the mesh $\mesh$. It is composed of $\ncell$ cells $K$ 
which in our implementation consists in tetrahedra.
There are $\nface$ faces $\face$ in the mesh, which are separated 
into the interior ones (interfaces between two cells) $\setface^i$, 
and the boundary ones $\setface^b$ on $\ddomain$ such that, 
\begin{equation}\label{eq:mesh-and-face}
  \mesh = \bigcup_{e=1}^{\ncell} K_e \,, \qquad \qquad 
  \setface = \bigcup_{e=1}^{\ncell} \partial K_e  = \setface^i \cup \setface^b \,=\, \bigcup_{k=1}^{\nface} \face_k \,.
\end{equation}

\paragraph{Jump}

The jump of a vector $\bv$ between two adjacent 
cells $K^+$ and $K^-$ is denoted by $\llbracket\, \cdot\, \rrbracket$ such that
\begin{equation}\label{eq:jump}
  \llbracket\, \bv \cdot \n \, \rrbracket \,=\,
  \bv^+ \cdot \n^+ \,+\, \bv^- \cdot \n^-,
\end{equation}
where $\n^{\pm}$ denotes the outward pointing normal,
and $\bv^{\pm}$ are the interior Dirichlet traces of 
$\bv|_{K^{\pm}}$ on $\partial K^{\pm}$.

\paragraph{Indexes}
Following \cref{eq:mesh-and-face}, a cell $K_e$ of the mesh uses index $e \in \{1,\ldots,\ncell \}$. 
We introduce two ways to index a face, the global index $\face_k$, for $k \in \{1,\ldots,\nface\}$, 
and the \emph{local} one $\face_{(e,\ell)}$ which means the face $\ell$ of element $e$:
\begin{equation}\begin{aligned}
&   \face_k \, ,  \hspace*{3em}  k \in \{1,\ldots,\nface\}, \hspace*{10em} \text{global indexing} \,; \\
&  \face_{(e,\ell)} \, , \qquad  e \in \{ 1,\ldots,\ncell\}\, , \quad
                               \ell\in  \{ 1,\ldots,\ncell^e\} \, \qquad \text{local indexing}\,.
\end{aligned} \end{equation}
Here, $\ncell^e$ is the number of faces for an element, i.e., four for tetrahedra.

\paragraph{Discretization spaces}

The space of discretization uses piecewise polynomial functions 
of order less than or equation to $\mathfrak{p}$, more specifically
we employ the Lagrange basis $\mathbb{P}_\mathfrak{p}$. The 
following spaces are used,
\begin{equation}\label{eq:fem-spaces} 
\left\lbrace \qquad \begin{aligned}
 W_h &= \left\lbrace  w_h \in L^2(\domain): 
        \hspace*{3.0em} w_h\!\mid_{K_e} \in \mathbb{P}_{\mathfrak{p}_e}(K_e), 
        \hspace*{2.4em} \forall K_e \in \mesh  \right\rbrace ;\\[0.3em]
\boldsymbol{W}_h= [W_h]^3 &= \left\lbrace 
      \bw_h \in L^2(\domain)^3 : \hspace*{1.50em} 
      \hspace*{1em}\bw_h \!\mid_{K_e} \in [\mathbb{P}_{\mathfrak{p}_e}(K_e)]^3, 
      \hspace*{1.2em} 
      \forall K_e \in \mesh  \right\rbrace; \\
V_h &= \left\lbrace  v_h \in L^2(\setface)\,:
                     \hspace*{3em} 
                     v_h \mid_{\face_k} \in \mathbb{P}_{\tilde{\mathfrak{p}}_k}(\face_k), 
                     \hspace*{3.4em} \forall \face_k \in \setface \right\rbrace.
\end{aligned} \right. 
\end{equation}
Different orders of approximation $\mathfrak{p}_e$ are allowed between each cell $K_e$, 
i.e., we use \emph{$p$-adaptivity}. In our implementation, the polynomial order 
$\tilde{\mathfrak{p}}_k$ on the face $k$ between elements  $K^+$ and $K^-$ is taken
as the maximum between the orders $\mathfrak{p}_+$ and $\mathfrak{p}_-$.

\subsection{HDG formulation} 

For the system \cref{eq:main-generic}, the HDG formulation 
is to find $(\bu_h, w_h, \lambda_h) \in \boldsymbol{W}_h \times W_h\times V_h $ such that: 
\begin{subequations}\label{eq:strong-form_volume}
\begin{empheq}[left={\text{On each cell $K_e \in \mathcal{T}_h$,}\,\,\,\empheqlbrace}]{align}
  \boldsymbol{A} \, \bu_h
  \,+\, \boldsymbol{\beta} \, w_h  
  \,+\, \nabla w_h \,=\, 0,  \quad &\text{ on } K_e\\
        \nabla\cdot \bu_h
    \,-\, \boldsymbol{\beta} \,\cdot\,\bu_h \,+\,\varrho\, w_h \,=\, h, \quad & \text{ on } K_e \\
    w_h = \lambda _h , \quad & \text{ on } \partial K_e.
\end{empheq}\end{subequations}
\begin{subequations}
\label{eq:strong-form_face}
\begin{align}
& \text{On each interior face \, $\face_k \, \in \, \setface^i$:} \qquad\qquad
  \llbracket  \yhwidehat{  \bu_h \cdot \n_{\face_k}} \rrbracket = 0 \,; \\[.5em]
& \text{On each boundary face \, $\face_k \, \in \, \setface^b$:} \qquad\qquad
  \dfrac{1}{\mathbf{Z}_{\mathrm{bc}} \cdot \n  }\lambda_h\, + \, \yhwidehat{ \bu_h\cdot \n} \,=\, 0\,,
\end{align}\end{subequations}
where we denote with $\widehat{\cdot}$ the numerical trace. 
The HDG numerical traces are further expressed in 
terms of a stabilization term $\tau$ such that:
\begin{equation}\label{eq:strong-form_jump}
\text{On each face \, $\face_k \, \in \, \setface$:} \qquad\qquad
\yhwidehat{ \bu_h \cdot \n_{\face_k}}  \,=\,  
            \bu_h\cdot  \n_{\face_k} \,-\, \tau \Big( w_h - \lambda_h \Big) \,.
\end{equation}

\begin{figure}[ht!]\centering
\includegraphics[scale=1]{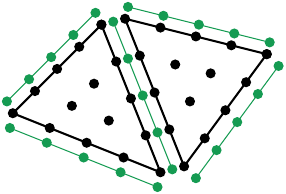}
\caption{Illustration of the degrees of freedom (dofs) for HDG discretization 
         in two dimensions with polynomial order $\mathfrak{p}=4$. The volume 
         unknowns are represented with the black dofs (unknowns $\bu_h$ and $w_h$, 
         hence the black dofs have multiplicity $(\mathrm{dimension}+1)$). 
         The green dofs (of multiplicity 1) are for the numerical 
         trace $\lambda$. 
         Only the dofs of the numerical traces are used in the global 
         linear system.}
\label{figure:hdg-dof-2D}
\end{figure}

In \cref{figure:hdg-dof-2D}, we illustrate the degrees of freedom 
for the HDG discretization in two dimensions to ease the visualization.
The volume unknowns, $\bu_h$ and $w_h$, are represented with the black
degrees of freedom and are independent per cell. 
The numerical trace $\lambda$ (green degrees of freedom) is instead 
only defined on the faces of the cells.

\paragraph{Weak form of the HDG formulation} 

To obtain the weak form, we integrate 
the volume problem \cref{eq:strong-form_volume} 
against test functions $(\bpsi,\varphi) \in \boldsymbol{W}_h\times W_h$,
and the problem on the faces 
\cref{eq:strong-form_face}
against test function $\zeta\in V_h$.
Using integration by parts and the definition of the 
HDG numerical trace \cref{eq:strong-form_jump}, we 
obtain the following problem, cf. \cite[Section~5]{Pham2025axisymreport} for the details.
\medskip

Find $(\bu_h, w_h,\lambda_h) \in \boldsymbol{W}_h \times W_h\times V_h$  such that,
\begin{itemize}
 \item On each cell $K_e \in \mathcal{T}_h$, 
       for all test functions $(\bpsi,\varphi) \in \boldsymbol{W}_h \times W_h$,
\begin{subequations}\label{eq:weak-form_volume}
\begin{empheq}[left={\empheqlbrace}]{align}
&       \int_{K}  \boldsymbol{A} \, \bu_h \, \bpsi \, \dd \bx
  + \int_{K}  w_h \boldsymbol{\beta}_1 \cdot \bpsi \,\dd \bx
  + \int_{\partial K} \lambda_h \, \bpsi \cdot \n  \,\dd \bx
  - \int_{K} w_h \nabla\cdot \bpsi \,\dd \bx \,=\, 0 \\[.50em]
&   \int_{K} \nabla\cdot \bu_h \varphi \, \dd \bx 
    \,-\,\int_{\partial K} \tau (w_h - \lambda_h) \varphi \, \dd \bx
    \,+\, 
    \int_{K} \boldsymbol{\beta}_2 \,\cdot\,\bu_h \varphi \, \dd \bx 
    \,+\,\int_{K} \varrho\, w_h \varphi \, \dd \bx \\
& \hspace*{25em} \, = \, \int_{K} f \, \varphi \,\dd \bx \,. \nonumber
\end{empheq}\end{subequations}
 \item On each interior face $\face_k \in \setface^i$, 
       for all test functions $\zeta \in V_h$:
\begin{equation}\label{eq:weak-form-face-interior}
  \int_{\face} \left( \bu_h^+ \cdot \n_\face^+  -
 \tau^+ (w_h^+ - \lambda_h)  \right) \zeta \,\dd\bx
+ \int_{\face} \left(\bu_h^- \cdot \n_\face^-  -
     \tau^- (w_h^- - \lambda_h) \right) \zeta\,\dd\bx = 0.
\end{equation}
 \item On each boundary face $\face_k \in \setface^b$, 
       for all test functions $\zeta \in V_h$:
\begin{equation}\label{eq:weak-form-face-exterior}
  \int_{\face} \left( \bu_h\cdot \n_\face  
             - \tau w_h  + \left(\tau + \dfrac{1}{\mathbf{Z}_{\mathrm{bc}}\cdot \n}\right) 
               \lambda_h \right)\,\zeta \, \dd \bx \,=\, 0 \,.
\end{equation}
\end{itemize}

\paragraph{HDG Stabilization}

It is essential to select the stabilization parameter $\tau$ efficiently 
in the HDG formulation, cf. \cite{Pham2024stabilization} and references 
therein. In this work, we follow the approach we introduced 
in \cite{Pham2025axisymreport}, where the same wave equation is 
considered in an axially symmetric setting, and where different 
choices of $\tau$ are investigated numerically. The current 3D problem 
builds upon those findings and for the numerical trace relation 
\cref{eq:strong-form_jump}, we use $\tau$ such that, 
\cite[Sections 5.5 and 6.3]{Pham2025axisymreport},
\begin{equation}
  \tau \, := \, \left\vert \Amat^{-1} \boldsymbol{\beta} \cdot \n \right\vert 
       \,-\, \num{e6}\,\ii \,\omega\left\vert \n^\mathrm{t} \, \Amat^{-1} \, \n \right\vert\,.
\end{equation}

\subsection{Matrix descriptions of the problem} 

We denote the coefficients of the unknowns $\left(\bu_h:=(u_{xh},u_{yh},u_{zh}),\, w_h\right)$ 
in the Lagrange basis functions $\varphi^e_j$ by $\mathrm{p}_j^{[\bullet]}$, 
with $\bullet$ indicating the corresponding unknown:
\begin{equation}
  w_h\!\mid_{K_e} \,=\, \sum_{j=1}^{\ndof^e} \mathrm{p}^{[w]}_j \, \varphi_j^e \,; \qquad 
  u_{x h}\!\mid_{K_e} \,=\, \sum_{j=1}^{\ndof^e} \mathrm{p}_j^{[u_x]} \, \varphi_j^e \,; \qquad 
  \text{similarly for } u_{y h}\!\mid_{K} \text{ and } u_{z h}\!\mid_{K}.
\end{equation}
Here, the polynomial order is allowed to vary with the cells, 
and we denote with $\ndof^e$ the corresponding number of degrees
of freedom on cell $K_e$.
The coefficients $\mathrm{p}_j^{[\bullet]}$ associated with a 
cell $K_e$ are represented into the vector $\mathsf{U}^e$ such
that,
\begin{equation}
\mathbf{U}^e \,  =\, 
\begin{pmatrix}  \mathbb{p}^e_{[u_x]}  &  \mathbb{p}^e_{[u_y]} &   \mathbb{p}^e_{[u_z]} &  \mathbb{p}^e_{[w]} \end{pmatrix}^t, \hspace*{4em} 
\text{ with } \hspace*{0.1cm}   \mathbb{p}^e_{[\bullet]} 
\, =\, \left(  \mathrm{p}^{[\bullet]}_j \right)_{j=1}^{\ndof^e}. 
\end{equation}

The HDG trace unknown $\lambda_h$
is represented on each face $\face_k \in \setface$ in the 
basis functions $\zeta_j^k$ with coefficients $\mathrm{q}_j^k$ such that,
\begin{equation}
 \lambda_h\!\mid_{\face_k} \,=\, \sum_{j=1}^{\ndofface}\, \mathrm{q}_j^k \, \zeta^k_j, \quad \text{with }
 \mathrm{q}_j^k = \mathrm{q}_j^{(e,\ell)} \hspace*{0.1cm} \text{ where }  \hspace*{0.1cm}\face_k = \face^{(e,\ell)}.
\end{equation}
We gather the coefficients $\mathrm{q}_j^k$ into the global vector 
$\Lambda$ of length $N_{\mathrm{dof}}^{\mathrm{face}}$ which is the 
total number of degrees of freedom considering all faces. We denote,
\begin{equation}
\Lambda = \left( \mathbb{q}^k\right)_{k=1}^{\nface} , \hspace*{0.2cm} \text{ with } \hspace*{0.2cm}  \mathbb{q}^k = ( \mathrm{q}_j^k )_{j=1}^{\ndofface}.
\end{equation}
Its restriction to an element $K_e$ is given by the matrix $\mathbb{R}_e$:
\begin{equation}
   \Lambda^e := \hdgR_e \Lambda\,=\, \begin{pmatrix}\mathbb{q}^{(e,\ell)} \end{pmatrix}_{\ell=1}^{\nface}, 
   \quad \text{ with } \quad \mathbb{q}^{(e,\ell)}  = \left( \mathrm{q}_j^{(e,\ell)}\right)_{j=1}^{\ndofface}.
\end{equation}
Conversely, we note that $\hdgR_e^t$ assigns the 
contribution defined on an element $K_e$ to the 
global vector of unknown $\Lambda$, cf. \cite[Equation (3.39)]{Pham2024stabilization}.

The local problem  \cref{eq:weak-form_volume} is written 
in matrix form such as, for each $K_e$, $1\leq e\leq \ncell$:
\begin{equation}\label{eq:hdg-matrix-system-volume}
  \hdgA^e \, \mathbf{U}^e \,+\, \hdgC^e \, \hdgR_e \Lambda \,=\, \hdgF^e \,.
\end{equation}
The discretization of the edge problems 
\cref{eq:weak-form-face-interior,eq:weak-form-face-exterior}
on $\setface$ is written as,
\begin{equation} \label{eq:hdg-matrix-system-interface}
\sum_{e=1}^{\ncell} \hdgR_e^t \left( \hdgB^e \, \mathsf{U}^e \,+\, \hdgL^e \, \hdgR_e \, \Lambda\right) \, = \,  0.
\end{equation}
The notation follows \cite{Faucher2020adjoint,Pham2024stabilization,Pham2025axisymreport}
and a description of the matrices can be found in \cite{Pham2025axisymreport} for 
axially symmetric background models.
Upon assuming that the matrices $\hdgA^e$ are invertible, 
one can express $\mathbf{U}^e$ in terms of $\Lambda$ with
\cref{eq:hdg-matrix-system-volume}:
\begin{equation}\label{eq:hdg-U} 
\mathbf{U}^e \, =\,   \left(\hdgA^e\right)^{-1} 
                      \left( \hdgF^e \, -\, \hdgC^e\, \hdgR_e \, \Lambda \right) 
                      \,, \hspace*{2em} \forall\, e =1, \ldots, \ncell.
\end{equation}
Replacing in \cref{eq:hdg-matrix-system-interface}, we obtain the 
(sparse) global problem in terms of $\Lambda$ only:
\begin{equation} \label{eq:hdg-global-linear-system}  
\mathbb{K} \,  \Lambda \, = \, \hdgS 
  \hspace*{2em}\text{where}\, \left\{\begin{array}{l}
  \mathbb{K}:=\hspace*{0.3cm} \displaystyle \sum_{e=1}^{\ncell}\,\hdgR^t_e\,\left( 
              \hdgL^e\, -\,\hdgB^e(\hdgA^e)^{-1}\, \hdgC^e\right)\,  \hdgR_e \,,\\[2em]
  \hdgS := -\displaystyle \sum_{e=1}^{\ncell} \hdgR^t_e   \,   \mathbb{B}^e \, (\hdgA^e)^{-1} \, \hdgF^e \,.
\end{array}\right.
\end{equation}
Therefore, the global linear system \cref{eq:hdg-global-linear-system} 
has the size of the number of degrees of freedom of $\Lambda$.

\section{Implementation for 3D helioseismology} 
\label{section:implementation}

In this section, we describe the specific aspects of the 3D wave problem 
in helioseismology and the resulting large-scale matrices that must be handled.
In \cref{subsection:strategy-mesh}, we present the profiles of the 
solar background model S \cite{christensen1996current,AtmoI2020}, along 
with our workflow for generating the 3D discretization mesh.
\cref{subsection:padaptivity} highlights the use of $p$-adaptivity 
to control the number of unknowns based on the local wavelength relative 
to the cell size. Next, we report the size of the resulting HDG global 
matrices for different frequencies and mesh configurations.
The HDG method is implemented in the open-source 
software \haven~(\url{https://ffaucher.gitlab.io/hawen-website/}), 
\cite{Hawen2021}, using \texttt{MPI} and \texttt{OpenMP} 
parallelisms. The contributions to the HDG local matrices 
(such as $\hdgA^e$, $\hdgC^e$ in \cref{eq:hdg-matrix-system-volume}) are independently created on each 
cell, allowing for massively parallel assembly of the global matrix.
The resulting global sparse linear system \cref{eq:hdg-global-linear-system}, 
is solved with the direct solver \mumps~\cite{Amestoy2001,Mumps2024} 
(\url{https://mumps-solver.org/index.php}), discussed \cref{section:mumps}.
Once the global system is solved, the solution vector $\Lambda$ 
is used to reconstruct the volume unknown $\mathbf{U}$ via \cref{eq:hdg-U}. 
This step involves solving on each mesh cell a dense linear system 
which is however small in size and completely independent across cells. 
Thus, this step is embarrassingly parallel and we use the 
library LAPACK, \cite{LAPACK} in our implementation.

In the numerical experiments, two supercomputers have been used: 
the cluster \adastra~at \texttt{CINES} (\url{https://www.cines.fr/calcul/adastra/})
and the cluster \mps~at the Max Planck Institute in G\"ottingen (MPS).
The cluster \adastra~is used for the most demanding experiments
as it has much more nodes, each of them composed of 2 AMD Genoa EPYC 
9654 with 96 cores, 2.4 GHz processors (3.7 GHz boost) with 768 Gio 
of DDR5-4800 MHz memory per node. 
The \mps~cluster is of smaller size, but allows to run experiments 
with the coarsest mesh, and each node is equipped with 
2 AMD EPIC 7763 with 64 cores and 1TB RAM DDR4 3200 MT/s per node,
hence more memory per node compared to \adastra.

\subsection{Solar background and mesh creation}
\label{subsection:strategy-mesh}

\begin{figure}[ht!]\centering
\includegraphics[scale=1]{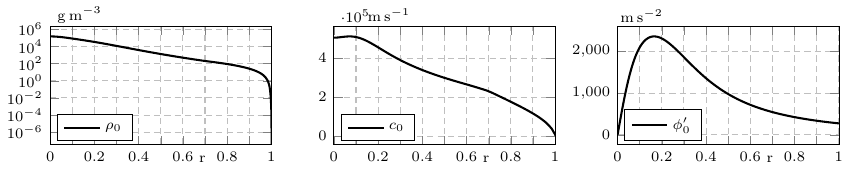}
\caption{The standard solar backgrounds from model S are radial 
\cite{christensen1996current}, i.e., they only depend on the distance 
to the origin $r=\Vert\bx\Vert$. The models are provided on a scaled
interval such that $r=1$ corresponds to the surface. The computations
are performed up until the maximal radius $\rmax=1.001$, \cite{Pham2025axisymreport}.}
\label{figure:modelS}
\end{figure}

We use the standard solar background model S \cite{christensen1996current,AtmoI2020} 
and consider a computational domain that extends slightly above the solar 
surface, as observations are typically made in the photosphere.
For this purpose, we employ a spherical domain with a scaled 
radius $\rmax = \num{1.001}$ (corresponding to about \num{696}\si{\kilo\meter} 
in the solar atmosphere), consistent with configurations used 
in solar modeling studies assuming 
spherical symmetry \cite{Pham2020Siam,Pham2024assembling} 
and axial symmetry \cite{Pham2025axisymreport}.
Note that we use a length-scale factor which corresponds 
to the Sun's radius (\num{696e6}\si{\meter}), and therefore 
$r=1$ corresponds here to the solar surface, cf.~\cite{Pham2020Siam,AtmoI2020,Pham2021Galbrun,Pham2024assembling,Pham2025axisymreport}.
The solar radial backgrounds are shown in \cref{figure:modelS}, 
with the density $\rho_0$, wave speed $c_0$ and the 
gravitational potential $\phi_0$. 
The challenge associated with the solar background
is the drastic decrease of the physical parameters 
near the solar surface, in particular
the density that decreases exponentially in the 
solar atmosphere.
This motivates the use of the variant equations that 
we introduced in \cref{section:wave-eq} 
(\cite{Pham2020Siam,Pham2024assembling,Pham2025axisymreport}), 
in which the density is not explicitly used, 
but its inverse scale-height (which is constant in 
the atmosphere) instead. 

\paragraph{Pipeline for mesh creation} 
To create the mesh, we rely on the open-source 
tools \texttt{TetGen} \cite{Hang2015Tetgen} 
and \texttt{GMSH} \cite{Geuzaine2009GMSH}. 
The discretization of the sphere with radius $\rmax = \num{1.001}$
is carried out with tetrahedral elements. Due to the decrease in the 
background model parameters near the surface, hence of the wavelength,
smaller elements are employed in the outer layers
(although $p$-adaptivity offers relative flexibility, see 
\cref{subsection:padaptivity}).
The workflow for the mesh creation is as follows:
\vspace*{-0.75em}

\begin{enumerate} \setlength{\itemsep}{-1pt}
  \item[1.] We list interior sphere radii 
            $r_{\mathrm{in},k}$ for $k=1,\ldots,n_1$, 
            with $r_{\mathrm{in},n_1}=\num{0.97}$, and $r_{\mathrm{in},k}$
            calculated from the radial solar 
            wavelength, \cite{Pham2021Galbrun}. 
            The list is then used to generate a concentric 
            tetrahedral mesh with \texttt{GMSH}. 
  \item[2.] From the mesh created in step 1, we retrieve 
            the current boundary faces (i.e., at height $r_{\mathrm{in},n_1}=\num{0.97}$), 
            and project at higher heights $r_{\mathrm{surf},k}$ for $k=1,\ldots,n_2$ 
            with $r_{\mathrm{surf},n_2} = \num{1.001}$.
  \item[3.] The final mesh is generated combining the interior 
            concentric mesh of step 1 and imposing the additional
            face layers generated in step 2, using \texttt{TetGen}.
\end{enumerate}

The step 2 consists in enforcing narrow layers where the spherical 
background models vary the most. 
These small layers cannot be directly implemented in step 1
because we have observed leakage between the layers: the surfaces 
are not respected in the triangulation when the radii are too close 
from one to the next. Imposing the faces as an input to \texttt{TetGen}
solves this issue. In the numerical experiments, we consider three meshes,
that only differ in step 2 (i.e., in the number of layers between 
\num{0.97} and \num{1.001}):
\vspace*{-0.75em}

\begin{itemize}\setlength{\itemsep}{-2pt}
  \item \meshA: mesh with  \num{6265192} elements;
  \item \meshB: mesh with  \num{8518333} elements;
  \item \meshC: mesh with \num{11525815} elements.
\end{itemize}
\begin{remark}
The mesh resulting from step 1. (i.e., up until the height $r=0.97$) 
has about 2 millions cells, which means that even in the smallest
mesh, \meshA, already two thirds of the elements are concentrated in
the interval $(0.97,1.001)$. 
\end{remark}

\subsection{$hp$-adaptivity and dimensions of global matrices}
\label{subsection:padaptivity}

\begin{figure}[ht!] \centering
  \includegraphics[scale=1]{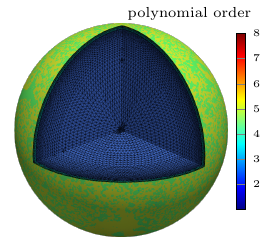} \hspace*{1em}
  \includegraphics[scale=1]{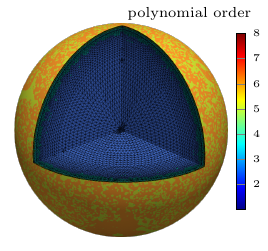}\hspace*{1em}
  \includegraphics[scale=1]{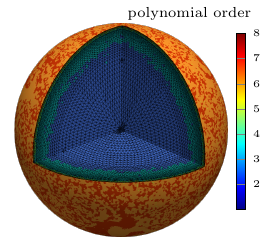}

  \caption{$hp$-adaptivity using \meshA~and frequencies 
           \num{1.5}\si{\milli\Hz} (\textbf{left}),
           \num{2}\si{\milli\Hz} (\textbf{middle}),
           and \num{2.5}\si{\milli\Hz} (\textbf{right}).
           The cells are refined near surface, and the 
           polynomial orders varying between 2 and 8 depending 
           on the local-to-the-cell wavelength. For visualization, 
           only the interior edges are shown.
           }
  \label{figure:mesh-p-adaptivity}
\end{figure}

Discretization methods in the Discontinuous Galerkin family are particularly 
adequate for $hp$-adaptivity, which means that each cell of the discretization
can have a different size and its own order of polynomial. As indicated above,
for $h$-adaptivity, we use relatively coarse cells in the interior, and refined
ones near surface. Then, the orders of the polynomials are determined 
based on the local-to-the-cell wavelength, in order to ensure a sufficient number 
of degrees of freedom (dofs). 
Empirically, we impose at least 10 dofs per wavelength for polynomial 
order 3, and at least 6 dofs per wavelength at order 8.
The $hp$-adaptivity is illustrated in \cref{figure:mesh-p-adaptivity} with \meshA~for
frequencies \num{1.5}, \num{2} and \num{2.5}\si{\milli\Hz}.
We observe that despite the mesh refinement near the surface, this region still 
shows the highest polynomial orders, whereas cells in the deeper interior can 
remain at order 2 or 3, even at 2\si{\milli\Hz} frequency.

\begin{remark}[Model representation on the mesh] \label{remark:model-representation}
  Employing high-order polynomials allows the use of relatively large cells 
  (in terms of local wavelength), which necessitates careful representation 
  of the background models to avoid having inaccurate representations, 
  \cite[Figure 7]{Faucher2020adjoint}. 
  For example, a piecewise-constant approximation per cell would be insufficiently 
  accurate in the near-surface layers of the Sun, where the background parameters 
  vary rapidly. 
  To address this, we represent the background models within each cell using 
  Lagrange basis. 
  Note that by eliminating the explicit dependence on density, the Liouville 
  variants of \cref{subsection:liouville-equations} also avoid the need to 
  project its exponentially decreasing profile onto the computational mesh.
\end{remark}

\paragraph{Size of the HDG global matrix} 

\begin{figure}[ht!] \centering
  \includegraphics[scale=1]{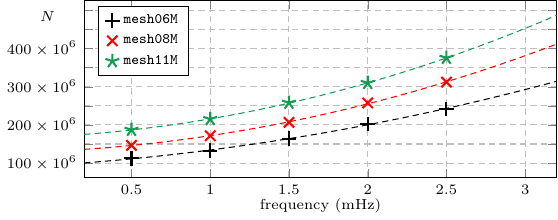}
  \vspace*{-0.75em}
  \caption{Size of the HDG global matrix with frequency and mesh
          (the numbers are given \cref{table:matrix-analytics}). The dashed lines
          correspond to the best quadratic approximation function.}
  \label{figure:hdg-matrix-size-function}
\end{figure}

In \cref{table:matrix-analytics,figure:hdg-matrix-size-function}, 
we provide the dimensions
of the HDG global matrix depending on the different meshes
and the frequency of the simulations. We use the 
$p$-adaptivity strategy in each case, i.e., each cell has
its polynomial order determined from the local wavelength,
and the polynomial order varies between 2 and 8.
As the frequency increases, the wavelength decreases, necessitating 
higher polynomial orders and consequently leading to larger problem sizes.

\begin{table}[ht!] \begin{center}
\caption{HDG problem dimensions depending on the mesh and frequency, 
         the polynomial orders are automatically selected depending on the 
         local-to-the-cell wavelength, and vary between \num{2} and \num{8},
         as illustrated in \cref{figure:mesh-p-adaptivity}.}
\vspace*{-0.50em} \scriptsize
\label{table:matrix-analytics}
\renewcommand{\arraystretch}{1.2}
\begin{tabular}{|>{\arraybackslash}p{.12\linewidth}|
                 >{\arraybackslash}p{.10\linewidth}|
                 >{\arraybackslash}p{.10\linewidth}|
                 >{\arraybackslash}p{.10\linewidth}|
                 >{\arraybackslash}p{.10\linewidth}|
                 >{\arraybackslash}p{.10\linewidth}|
               }

\hline
\multicolumn{6}{|c|}{Size $N$ of the HDG global matrix $\hdgK$ \cref{eq:hdg-global-linear-system}} \\ \hline
  & \multicolumn{1}{ c|}{\num{0.5}\si{\milli\Hz}}
  & \multicolumn{1}{ c|}{\num{1}\si{\milli\Hz}} 
  & \multicolumn{1}{ c|}{\num{1.5}\si{\milli\Hz}}
  & \multicolumn{1}{ c|}{\num{2}\si{\milli\Hz}} 
  & \multicolumn{1}{ c|}{\num{2.5}\si{\milli\Hz}}\\ \hline
\meshA   & 
\num{112088015} &
\num{132231488} &
\num{163806169} & 
\num{202431918} &
\num{242128419} 
\\ \hline
\meshB   & 
\num{146559348} &
\num{171584711} &
\num{206630506} &
\num{257576619} &
\num{312100698} 
\\ \hline
\meshC          & 
\num{187022492} &
\num{215289472} &
\num{257953840} &
\num{309859823} &
\num{375747408} 
\\ \hline
\multicolumn{6}{c}{ } \\ \hline

\multicolumn{6}{|c|}{Size of the HDG global matrix 
                     compared to the total number 
                     of volume unknowns,} \\
\multicolumn{6}{|c|}{i.e.,
                     $\mathrm{dim}(\Lambda)/(\sum_e\mathrm{dim}(\mathbf{U}^e))$.} \\ \hline
  & \multicolumn{1}{ c|}{\num{0.5}\si{\milli\Hz}}
  & \multicolumn{1}{ c|}{\num{1}\si{\milli\Hz}} 
  & \multicolumn{1}{ c|}{\num{1.5}\si{\milli\Hz}}
  & \multicolumn{1}{ c|}{\num{2}\si{\milli\Hz}} 
  & \multicolumn{1}{ c|}{\num{2.5}\si{\milli\Hz}} 
  \\ \hline
\meshA   &     
\num{25.68}\si{\percent} & 
\num{24.73}\si{\percent} & 
\num{22.95}\si{\percent} & 
\num{21.08}\si{\percent} & 
\num{19.51}\si{\percent} 
\\ \hline  
\meshB   &
\num{26.07}\si{\percent} & 
\num{24.90}\si{\percent} & 
\num{23.66}\si{\percent} & 
\num{21.79}\si{\percent} & 
\num{20.10}\si{\percent}  
\\ \hline
\meshC   &  
\num{26.85}\si{\percent} & 
\num{25.69}\si{\percent} & 
\num{24.26}\si{\percent} & 
\num{22.59}\si{\percent} & 
\num{20.99}\si{\percent}   
\\ \hline
\end{tabular}
\end{center}
\end{table}

The sizes of matrices vary from 100 to 400 millions, requiring for 
high-performance computing infrastructure for their handling. The coarsest
mesh \meshA~naturally leads to the smallest size of matrix, as it avoids 
over-refined areas (i.e., areas where the number of points per wavelength is 
high), and compensate coarse cells by employing higher-order polynomials.
In \cref{figure:hdg-matrix-size-function}, we see that for each mesh, 
a quadratic function captures the evolution of the matrix size 
accurately, that is, there is a $\omega^2$ dependence on the size 
of the matrix. This allows us to anticipate the size of the system 
to handle with frequencies. For instance, at \num{3}\si{\milli\Hz},
the matrix size should be of about \num{300} and \num{400} millions,
respectively for \meshA~and \meshC.

  In \cref{table:matrix-analytics}, 
  the threshold of 25\si{\percent} between the size of the global 
  matrix and number of volume unknowns indicates that the HDG 
  method leads to a significant reduction in the size of the linear 
  system compared to the standard DG methods. 
  Indeed, the original system of equations \cref{main-equation:original-div_v0} 
  can be reformulated as a second-order partial 
  differential equation involving only the Eulerian pressure $\pressE$,
  leading to an anisotropic Helmholtz equation with a scalar unknown
  (cf. \cite{Pham2025axisymreport})
  effectively decreasing the number of unknowns by a factor 
  of four by eliminating the vector-valued variable. 
  The static condensation inherent to the HDG method allows us to 
  retain the vector unknowns while working with a global linear system 
  in terms of one unknown ($\Lambda$) which, in addition, ignores the
  interior-to-the-cell dofs.
  In contrast, a DG discretization applied directly to the scalar equation would 
  require additional post-processing to recover these vector quantities. 
  
  Moreover, we see in \cref{table:matrix-analytics} that for the highest frequencies, 
  the size of  the HDG global matrix remains consistently below 25{\si{\percent}} 
  compared to the number of volume unknowns, i.e., the system is smaller than considering
  a scalar problem with DG discretization method. 
  Note that comparing these numbers with a continuous Galerkin (CG) discretization 
  is less straightforward due to the difficulties in implementing $p$-adaptivity in 
  CG methods, and the need to consider the solution accuracy for a fair comparison.

\section{Efficient linear system solution with \mumps}
\label{section:mumps}

The efficiency of our approach relies on the ability to 
efficiently solve the global sparse linear system resulting 
from the HDG discretization (of sizes given in 
\cref{table:matrix-analytics}) with multiple right-hand sides (RHS).
In this study, we use a direct method with the scientific library 
\mumps~\cite{Mumps2024}.
A primary motivation for choosing a direct solver is the need to efficiently handle 
multiple right-hand sides, as commonly encountered in helioseismology. 
Moreover, even though iterative solvers are often more memory-efficient, 
their performance relies on the availability of a robust and problem-specific 
preconditioner, whose design is challenging in the context of stellar models.
To control the complexity (in terms of number of operations and memory footprint) of the numerical phases of the solution,
block low-rank compression  with 
mixed-precision storage, referred to as BLR in the remainder of this paper, 
of the factor matrices \cite{aabblw:15,ablm:17b,Amestoy2023mixed} are used.
We analyse in this section the efficiency of the proposed approach for 3D helioseismology.

\subsection{Main features and performance indicators for large-scale simulations}
\label{subsection:mumps-features}

Considering the sparse linear system $AX = B$, where $B$ is a 
matrix of right-hand sides (RHS), which is a sparse matrix in 
our case, the solution $X$ is obtained with \mumps~in three main 
steps:
\begin{enumerate}\setlength{\itemsep}{-1pt}
  \item \textbf{Analysis phase}: the matrix $A$ is permuted to reduce 
        the amount of fill-in in the factors and a symbolic factorization 
        step is performed to determine all internal data structures.

  \item \textbf{Factorization phase}: the matrix $A$ is factorized as
        $A=LU$ (as $A$ is unsymmetric, otherwise one uses 
        $A=LDL^T$ for symmetric matrices). 

  \item  \textbf{Solve phase}: it computes the solution for each RHS via
         a forward elimination 
         $L Y = B$ 
         followed by a backward substitution 
         $U X = Y$ to compute the solution $X$.
         Note that the solution of multiple sparse 
         RHS is performed efficiently once the analysis 
         and factorization are completed.
\end{enumerate}

Firstly, the condition number of the matrices associated 
with the different variants of equations in 
\cref{eq:equation-system-notation} is analyzed in 
\cref{subsection:cond-number}. We then retain 
the formulation that yields the most well-conditioned systems.
Secondly, the performance of the three phases is evaluated measuring both the
computational time and the memory consumption, while controlling 
the accuracy of the solution.
In the following, we describe in more details four features of the linear solver that
strongly impact the performance of our simulation.

\begin{itemize}[leftmargin=*]
\item \textbf{Analysis using block format/Analysis by block.}
\mumps~exploits the specific structure of the global matrix 
resulting from the HDG discretization. 
The graph $G_A$ associated to the original matrix $A$ is defined 
such that each column/row of $A$ corresponds to a node in 
$G_A$, and each non-zero entry $a_{ij}$ corresponds to an 
edge between the nodes $i$ and $j$. 
As illustrated in \cref{figure:hdg-dof-2D}, the HDG dofs belonging 
to the same face share the same connectivity. It follows that all 
nodes associated with the same face have the same adjacency in $G_A$. 
During the analysis phase, $G_A$ can thus be replaced by a 
graph $G_{A\texttt{comp}}$ associated with the faces instead, 
thus reducing the number of nodes and edges in the graph. 
As a simple illustration, assuming each face has $k$ dofs, 
the number of nodes in $G_{A\texttt{comp}}$ is reduced by 
a factor $k$, and the number of edges by approximately 
a factor $k^2$.
In addition, the number of faces being independent 
of the polynomial order, the cost of the analysis phase 
is also independent of the frequency.

\item \textbf{Block low-rank.}
The complexity of the factorization is reduced through the use of 
block low-rank (BLR) approximation, developed by 
\mumps~group to exploit data sparsity \cite{aabblw:15,ablm:17b}.
It has been proved in \cite{ablm:17a} that the complexity 
of sparse solvers can be reduced in terms of both floating-point 
operations and memory usage, particularly during the factorization phase. 
In \cite{hima:21}, the BLR method was proven to be backward stable, 
meaning that the computed solution corresponds to the exact solution 
of a nearby problem with a small perturbation in the input data. 
This property is crucial, as it allows for
safely selecting a numerical threshold $\epsilonblr$ for accuracy.
In \cref{subsection:block-low-rank}, we use ${\texttt{blr-}}\epsilonblr$ 
to indicate the use of mixed precision and BLR with threshold $\epsilonblr$,
in particular comparing between \blr{5}, \blr{7} and \blr{9}.

Note that unlike hierarchical formats, the BLR format is based 
on a flat, non-hierarchical blocking of the matrix. In fact, 
we have shown in~\cite{ablm:17a} that despite being non-hierarchical, 
the BLR format can still achieve a substantial complexity reduction. Furthermore this flat format 
enables to implement all numerical features of a general sparse solver including numerical pivoting
to preserve backward stability. 

\item \textbf{Mixed-precision storage.}
Mixed precision has been introduced to reduce the memory footprint using up to seven precision to store the matrices of factors $L$ and $U$. In our context the matrices of factors are represented in BLR format at accuracy $\epsilonblr$. 
This low-rank approximation exhibits remarkable properties that enable the use of a continuum of precisions without losing accuracy. 
In our approach, we have decoupled the storage precision from the compute precision: the data are stored and accessed in low precision, but the computations are
performed in high precision.
Given a block $B$ of the $LU$ factors,  let us assume that for the clarity of our discussion that a singular value (SVD) decomposition of $B$ is performed. Dropping all singular values and singular vectors associated with singular values smaller than $\epsilonblr$ may lead to a more compact low-rank representation of $B$ at precision $\epsilonblr$. 
Furthermore, as demonstrated in \cite{Amestoy2023mixed}, 
each singular vector can be stored in a precision inversely proportional to its associated singular value while maintaining high accuracy overall.
This enables us to obtain a mixed-precision representation of the low-rank block $B$. Doing so on the BLR representation of the matrices of factors, one 
further reduces the memory footprint of the factorization.

In our mixed-precision algorithm for single precision (32-bits) factorization, three different precision formats (32-bit, 24-bit, and 16-bit) are used
used to store the $LU$ factors 
(whereas seven are used for double-precision factorization).
The mixed-precision representation of the $LU$ factors is computed 
during factorization. 
During the solution phase the mixed-precision representation of the 
factors is then converted to the computing precision (32-bits).


\end{itemize}


\paragraph{Indicators}
The indicators provided by \mumps~are employed to evaluate 
the computational cost of the linear system solution, 
offering insights into both time and memory consumption 
of the different phases. 
In particular, to assess the gain achieved by the BLR factorization 
and mixed-precision techniques, three indicators are used: 
\pctNop~and \pctNen~are respectively the percentage of the number 
of operations and entries to store during the BLR factorization 
compared to a full-rank (FR) factorization;  
\pctNMP~represents the gain in terms of the reduction of the storage 
of the $LU$ factors 
obtained with the mixed-precision representation during a BLR factorization
as compared to a BLR factorization without mixed precision.
Using mixed precision during BLR factorization does modify
the number of entries in the factors \pctNen.
\cref{table:mumps-analytics} lists the indicators, 
given in terms of the default outputs provided 
by \mumps \cite{mumpsuserguide}. 


\begin{table}[ht!] \begin{center}
\caption{List of indicators used to evaluate the performance of the linear 
         solver during the analysis and the factorization phases. Here we 
         respect the nomenclature of \mumps~for the notation of the 
         indicators, cf.~\cite{mumpsuserguide}, from which we define
         three main indicators: \pctNop, \pctNen~and~\pctNMP.}
\vspace*{-0.50em} \small
\label{table:mumps-analytics}
\renewcommand{\arraystretch}{1.1}
\begin{tabular}{|>{\arraybackslash}p{.22\linewidth}|
                 >{\arraybackslash}p{.73\linewidth}|
               }
\hline
\multicolumn{2}{|c|}{\textbf{Indicators for the analysis phase}} \\ \hline
Graph size      &  Originally, the graph has as many nodes as the number 
                   of columns/rows of the matrix. 
                   It is reduced when the block format is exploited.\\ \hline
Number of edges &  Originally, it is the number of non-zero entries in the matrix.
                   It is reduced when the block format is exploited. \\ \hline
Analysis time   &  Total time to perform the analysis phase. \\ \hline 
\end{tabular}
\vspace*{1em}

\begin{tabular}{|>{\arraybackslash}p{.22\linewidth}|
                 >{\arraybackslash}p{.73\linewidth}|
               }
\hline
\multicolumn{2}{|c|}{\textbf{Indicators for the factorization phase}} \\ \hline
Factorization memory &  The memory effectively used during factorization step. 
                        For robustness, the total memory allocated is in 
                        practice higher. It is given in Terabytes (TB). \\ \hline
Factorization time   & Total time to perform the factorization.  \\ \hline
\end{tabular}

\vspace*{0.85em}

\begin{tabular}{|>{\arraybackslash}p{.22\linewidth}|
                 >{\arraybackslash}p{.73\linewidth}|
               }
\hline
\texttt{RINFOG(3)}, \texttt{RINFOG(14)}  
& Number of floating point operations to perform the full-rank (\fr) 
  and BLR factorizations, respectively. 
                       \\ \hline
\pctNop$:=100\frac{\texttt{RINFOG(14)}}{\texttt{RINFOG(3)}}$ &  
                      Percentage of operations in the BLR factorization
                      compared to the FR factorization.
       \\ \hline
\end{tabular}

\vspace*{0.85em}

\begin{tabular}{|>{\arraybackslash}p{.22\linewidth}|
                 >{\arraybackslash}p{.73\linewidth}|
               }
\hline
\texttt{INFOG(29)}, \texttt{INFOG(35)}  
& Number of entries in factors for the FR and BLR factorizations, respectively. \\
\texttt{INFOG(9)}
& Number of entries used to store the factors. 
                       \\ \hline
\pctNen$:=100\frac{\texttt{INFOG(35)}}{\texttt{INFOG(29)}}$ &  
                      Percentage of entries in factors for the BLR factorization
                      compared to the FR factorization. \\ \hline
\end{tabular}

\vspace*{0.85em}

\begin{tabular}{|>{\arraybackslash}p{.22\linewidth}|
                 >{\arraybackslash}p{.73\linewidth}|
               }
\hline
\pctNMP$:=100\frac{\texttt{INFOG(9) }}{\texttt{INFOG(35)}}$ 
& Percentage of the storage of the factors with the BLR factorization with mixed precision
  compared to BLR factorization. \\ \hline
\end{tabular}
\end{center}
\end{table}

\subsection{Conditioning of the global matrix depending on the variants}
\label{subsection:cond-number}
\newcommand{\bwderr}{\mathfrak{b}_{\mathrm{wd}\epsilon}}

In \cref{section:wave-eq}, three variants of the system are introduced, 
respectively denoted as $\Loriginaldiv$, $\Lliouville$ and $\Lliouvillec$. 
Here we show the impact of the formulations on the global matrix, by 
using the conditioning indicators in \mumps \cite[Section 3.3.2]{mumpsuserguide}.
Two indicators are used to evaluate the behaviour of the system:
the condition number $\texttt{cond}$, and the componentwise backward error $\bwderr$ 
which takes also into account the sparsity of $b$ \cite{Arioli1989,mumpsuserguide}. 
Consider the general system $Ax=b$, with $b$ a RHS vector. 
Due to limited numerical precision and/or the use of the BLR 
factorization, the computed solution 
$\tilde{x}$ is the exact solution to an approximate problem,
\begin{equation}
   (A \,+\, \delta_A) \, \tilde{x} \,=\, (b\,+\, \delta_b) \,, \qquad\qquad
   \tilde{x} = x + \delta_{x} \,.
\end{equation}
The componentwise backward error $\bwderr$ is such that
\begin{equation}
  \dfrac{\vert \delta_{a_{ij}} \vert}{\vert a_{ij} \vert} \, \leq \, \bwderr \,,
\end{equation}
where $a_{ij}$ and $\delta_{a_{ij}}$ are components 
of $A$ and $\delta_A$, respectively. Note that $\bwderr$ 
is impacted by the BLR threshold used. 
Furthermore, an upper bound of the forward error is given by,
\begin{equation}
  \dfrac{\Vert \delta_{x} \Vert_{\infty}}{\Vert x \Vert_\infty} 
  \quad \leq \quad \texttt{cond} \,\, \times \,\, \bwderr \,.
\end{equation}
In \cref{table:cond}, we provide the values of $\texttt{cond}$ 
and $\bwderr$, comparing the different variants.

\begin{table}[ht!] \begin{center}
\caption{Condition numbers depending on the 
         variants \cref{eq:equation-system-notation} with frequency
         and mesh.
         The size of the matrices is given 
         in \cref{table:matrix-analytics}. 
         The experiments are carried out using \blr{7}.}
\vspace*{-0.50em} \scriptsize
\label{table:cond}

\sisetup{exponent-product = \hspace*{-0.10em}\times\hspace*{-0.10em}, 
round-mode = figures, 
round-precision = 3}

\renewcommand{\arraystretch}{1.10}

\begin{minipage}[t]{.49\linewidth}
\begin{tabular}{|>{\arraybackslash}p{.22\linewidth}|
                 >{\arraybackslash}p{.22\linewidth}|
                 >{\arraybackslash}p{.22\linewidth}|
                 >{\arraybackslash}p{.22\linewidth}|
               }
\hline
  & \multicolumn{3}{ c|}{frequency 0.5mHz -- \meshA} \\ \hline
  & $\Loriginaldiv$ &  $\Lliouville$ &  $\Lliouvillec$  \\ \hline
 $\texttt{cond}$           &\num{1.97D+08} & \num{1.05D+06} & \num{2.12D+06}  \\
 $\bwderr$                 &\num{6.39D-09} & \num{6.60D-09} & \num{5.93D-09}  \\ \hline
 $\texttt{cond}\times\bwderr$  &\num{1.25883}  & \num{0.00693}  & \num{0.0125716} \\\hline
\end{tabular}
\end{minipage} \hspace*{3em} \begin{minipage}[t]{.44\linewidth}
\begin{tabular}{|>{\arraybackslash}p{.22\linewidth}|
                 >{\arraybackslash}p{.22\linewidth}|
                 >{\arraybackslash}p{.22\linewidth}|
               }
\hline
    \multicolumn{3}{|c|}{frequency 1mHz -- \meshA} \\ \hline
    $\Loriginaldiv$ &  $\Lliouville$ &  $\Lliouvillec$  \\ \hline
\num{4.65D+09} & \num{2.61D+06} & \num{3.99D+06}  \\
\num{2.53D-08} & \num{2.33D-09} & \num{1.75D-09}  \\ \hline
\num{117.645}  & \num{0.0060813}& \num{0.0069825} \\ \hline
\end{tabular}
\end{minipage}

\vspace*{1em}

\begin{minipage}[t]{.49\linewidth}
\begin{tabular}{|>{\arraybackslash}p{.22\linewidth}|
                 >{\arraybackslash}p{.22\linewidth}|
                 >{\arraybackslash}p{.22\linewidth}|
                 >{\arraybackslash}p{.22\linewidth}|
               }
\hline
  & \multicolumn{3}{ c|}{frequency 1.5mHz -- \meshA} \\ \hline
  & $\Loriginaldiv$ &  $\Lliouville$ &  $\Lliouvillec$  \\ \hline
 $\texttt{cond}$          & \num{1.68D+10} & \num{7.16D+06} & \num{1.07D+07} \\ 
 $\bwderr$                & \num{2.14D-07} & \num{1.20D-08} & \num{2.69D-08} \\ \hline
 $\texttt{cond}\,\bwderr$ & \num{3595.2}   & \num{0.08592}  & \num{0.28783}  \\ \hline
\end{tabular}
\end{minipage} \hspace*{3em} \begin{minipage}[t]{.44\linewidth}
\begin{tabular}{|>{\arraybackslash}p{.22\linewidth}|
                 >{\arraybackslash}p{.22\linewidth}|
                 >{\arraybackslash}p{.22\linewidth}|
               }
\hline
    \multicolumn{3}{|c|}{frequency 1mHz -- \meshB} \\ \hline
    $\Loriginaldiv$ &  $\Lliouville$ &  $\Lliouvillec$  \\ \hline
\num{5.26e+09} & \num{2.63e+06} & \num{4.00e+06} \\ 
\num{1.42e-08} & \num{1.31e-08} & \num{1.00e-08} \\ \hline
\num{74.692  } & \num{0.034453} & \num{0.04}     \\ \hline
\end{tabular}

\end{minipage}

\end{center}
\end{table}

We observe the the condition number $\texttt{cond}$ 
for the original formulation $\Loriginaldiv$ is 
consistently at least two orders of magnitude higher
than that of the Liouville variants. 
Among the latter, $\Lliouville$ is  
slightly better than $\Lliouvillec$.
The backward error $\bwderr$  is relatively similar 
in all cases although the formulation $\Loriginaldiv$ 
is always slightly worse.
Consequently the forward error given by the 
product $\texttt{cond}\times\bwderr$ is the highest
for $\Loriginaldiv$ by several orders of magnitude, while 
$\Lliouville$ and $\Lliouvillec$ are in the same range, 
with $\Lliouville$ being slightly better. 
It means that, given an expected accuracy of the solution 
(represented by $\texttt{cond} \times \bwderr$), 
the Liouville variants are able to 
consider a larger $\epsilonblr$ (which will increase $\bwderr$) 
and thus more compression compared to $\Loriginaldiv$.
Additionally, the condition number tends to 
increase with frequency, as expected since the
wavelength shortens.
When comparing the different meshes, the condition number 
remains relatively close. 
We consequently select $\Lliouville$ for the subsequent experiments,
which is also more practical than $\Lliouvillec$ for modeling wave speed 
heterogeneities as the latter involves $\boldsymbol{\alpha}_c$ 
in \cref{main-equation:Liouville-c}, hence requiring derivatives 
of wave speed perturbations.


\subsection{Analysis by block}
\label{subsection:analysis-by-block}

We investigate the performance of the analysis phase, 
in particular using the analysis by block decribed in \cref{subsection:mumps-features}. 
In \cref{table:analysis-by-block}, we report the size 
of the graphs 
and the analysis times 
for different partitioners: \texttt{Scotch} \cite{scotch}
with and without analysis by block, and \texttt{Metis} \cite{metis}
with analysis by block. 
Originally (without exploiting the block format), the graph 
has the size of the matrix (\num{132231488}) and it 
has \num{11892113408} edges, requiring the use of 
\texttt{int64} for integer representation within the partitionner.
In contrast, exploiting the block format, the graph is 
reduced to size \num{12591470}, i.e., 
less than \num{10}\si{\percent} of the original size,
and fewer than \SI{1}{\percent} of the edges
remain (\num{75182304}).
As a result, enabling the analysis by block leads to a reduction 
in the analysis time of approximately \SI{70}{\percent}.
When comparing partitioners, we observe that \texttt{Scotch} is 
approximately twice as fast as \texttt{Metis} for our problem. 
For the subsequent experiments, we consequently exploit the block
format in the analysis, and employ the \texttt{Scotch} partitioner.

\begin{table}[ht!] \begin{center}
\caption{Information on the analysis phase for computations 
         uses \meshA~at frequency \num{1}\si{\milli\Hz}.
         }
\vspace*{-0.50em} \scriptsize
\label{table:analysis-by-block}

\renewcommand{\arraystretch}{1.20}

\begin{tabular}{|>{\arraybackslash}p{.20\linewidth}|
                 >{\arraybackslash}p{.20\linewidth}|
                 >{\arraybackslash}p{.20\linewidth}|
                 >{\arraybackslash}p{.20\linewidth}|
               }
\hline
  & Partitioner \texttt{Scotch} 
  & Partitioner \texttt{Scotch} 
  & Partitioner \texttt{Metis}  \\ 
  & without block-analysis 
  & with block-analysis 
  & with block-analysis \\ \hline 
Size of the graph     & \num{132231488} & \num{12591470} & \num{12591470}  \\ 
Number of graph edges & \num{11892113408} & \num{75182304} & \num{75182304} \\ 
Time for the analysis & \num{222}\si{\second} & \num{68}\si{\second} & \num{155}\si{\second} \\ \hline
\end{tabular}
\end{center}
\end{table}


\subsection{Block low-rank compression and mixed precision for factorization}
\label{subsection:block-low-rank}

The key features used in \mumps~to reduce the complexity 
(in terms of the number of operations and the memory footprint) of the 
matrix factorization are the block low-rank (BLR) compression and the 
mixed precision for the storage of the factors. 
The low-rank compression is applied with an input threshold 
$\epsilonblr$, which determines the expected backward error of the system to solve.
In this section, we explore the trade-off between memory consumption, 
which decreases when the BLR threshold increases, and the solution 
accuracy, which deteriorates as the threshold increases.
To evaluate the accuracy, we introduce the relative difference $\err_{xz}$ 
between solutions obtained from the full-rank (\fr) and BLR factorizations, 
evaluated on the $xz$-plane:
\begin{equation} \label{relative-difference:xz}
  \err_{xz}(\bx) 
  \,=\, \dfrac{\vert w_{\texttt{ref}}(\bx) - w_{\texttt{blr}}(\bx) \vert}{\Vert w_{\texttt{ref}}(\bx) \Vert} \,, 
  \hspace*{2em} \text{for $\bx$ in the $xz$-plane,} 
\end{equation}
Here the reference solution $w_{\texttt{ref}}$ corresponds to the 
one obtained with full-rank factorization.

\begin{figure}[ht!]\centering
          \includegraphics[height=2.70cm]{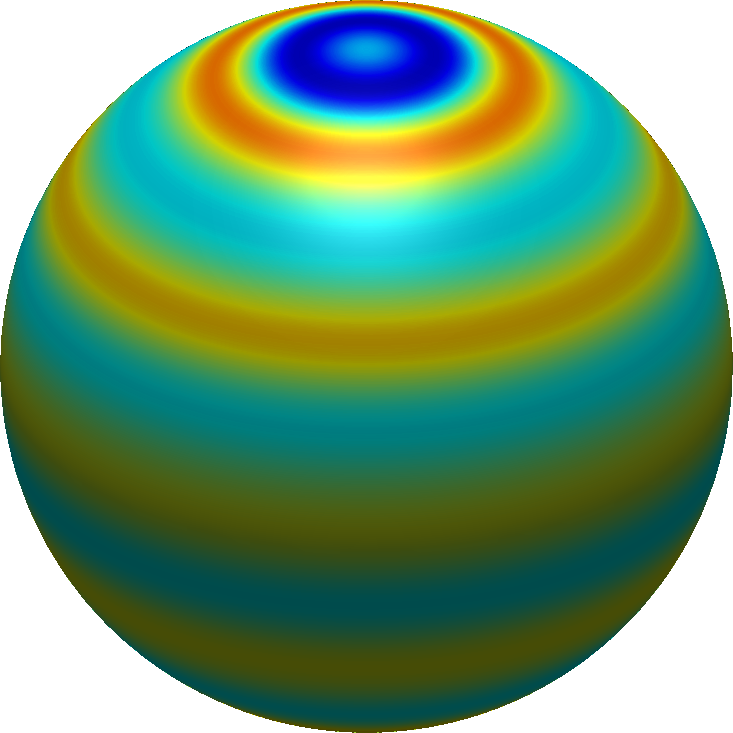}
          
  \pgfkeys{/pgf/fpu=true} 
  \pgfmathsetmacro{\cmin}{-7e-7} 
  \pgfmathsetmacro{\cmax}{7e-7} 
  \pgfkeys{/pgf/fpu=false}
  \setlength{\cbarwidth} {0.50em}
  \setlength{\cbarheight}{2.50cm}
  \hspace*{-4em}\raisebox{0.50em}{\input{figures/cbar}}

          \hspace*{-10em}\raisebox{-4.0em}{
          \begin{tikzpicture}[scale=1,every node/.style={font=\scriptsize,},
                              execute at begin node={\renewcommand{\baselinestretch}{1}}]
          \node[text width=3.21cm]{\scriptsize Wavefield at $r=1$ 
                                for a source at $(0,0,0.7)$ 
                                at frequency 1\si{\milli\Hz}.};
          \end{tikzpicture}}
          \includegraphics[scale=0.95]{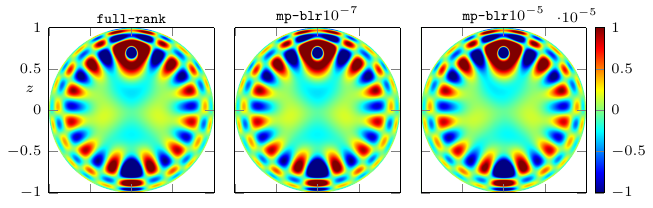} 
          \\[-4.75em]
\subfloat[Solutions $\mathrm{Re}(w)$ and relative difference for a source at $\zsrc=\num{0.7}$.]{
          \hspace*{3.85cm}\includegraphics[scale=0.95]{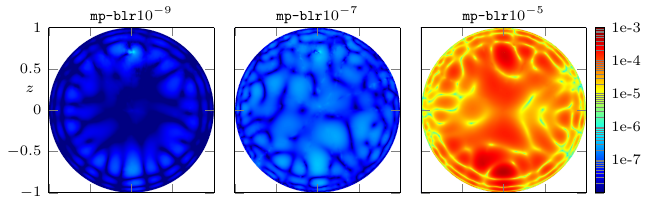}
          }

          \includegraphics[height=2.70cm]{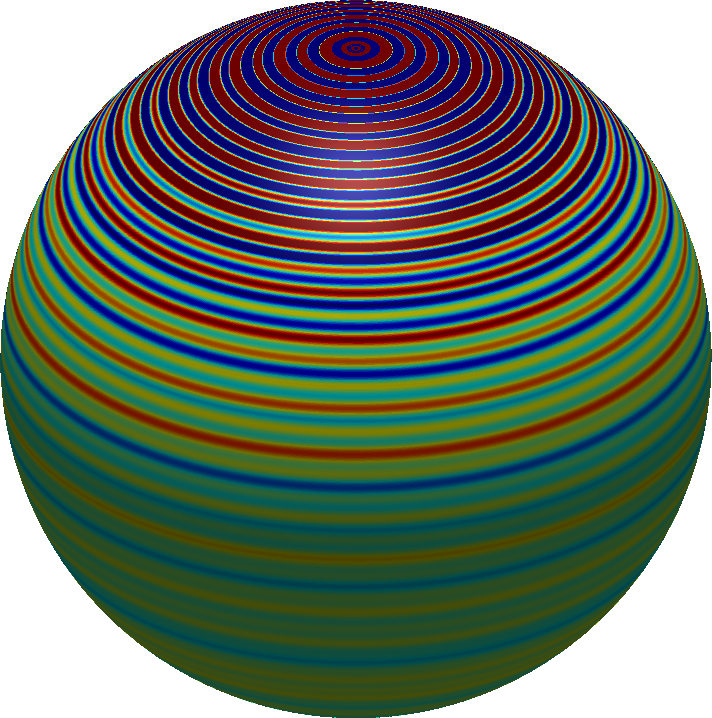}
          
  \pgfkeys{/pgf/fpu=true} 
  \pgfmathsetmacro{\cmin}{-5e-8} 
  \pgfmathsetmacro{\cmax}{5e-8} 
  \pgfkeys{/pgf/fpu=false}
  \setlength{\cbarwidth} {0.50em}
  \setlength{\cbarheight}{2.50cm}
  \hspace*{-4em}\raisebox{0.50em}{\input{figures/cbar}}

          \hspace*{-10em}\raisebox{-4.0em}{
          \begin{tikzpicture}[scale=1,every node/.style={font=\scriptsize,},
                              execute at begin node={\renewcommand{\baselinestretch}{1}}]
          \node[text width=3.21cm]{\scriptsize Wavefield at $r=1$ 
                                for a source at $(0,0,0.99986)$ 
                                at frequency 1\si{\milli\Hz}.};
          \end{tikzpicture}}
          \includegraphics[scale=0.95]{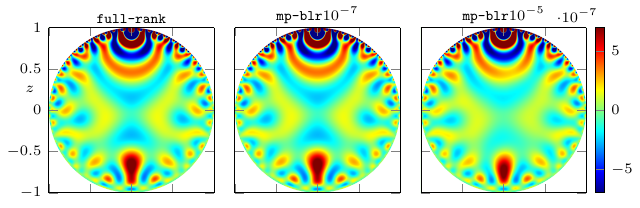}\\[-4.75em]
\subfloat[Solutions $\mathrm{Re}(w)$ and relative difference for a source at $\zsrc=\num{0.99986}$.]{
          \hspace*{3.85cm}\includegraphics[scale=0.95]{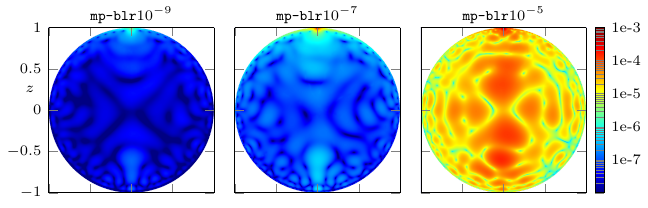}}
\caption{Simulations using \meshA~at frequency 1\si{\milli\Hz} for two different sources. 
         On the left, panels show the solution at a fixed radius $r=1$, and on the right panels are cuts on the $xz$-plane. On the right panels, the top row corresponds to the real part of the solution $w$ and
the bottom row to the relative error $\err_{xz}$ of \cref{relative-difference:xz}.
The experiments use the radial solar background model \cref{section:implementation} 
and attenuation $\gamma_{\mathrm{att}}/(2\pi)=\num{10}\si{\micro\Hz}$.}
\label{figure:wavefields-1mHz-fr-blr}
\end{figure}

In \cref{figure:wavefields-1mHz-fr-blr}, we compare the solutions 
for \meshA~at frequency 1\si{\milli\Hz} obtained with
\fr~and BLR factorizations with three different values of $\epsilonblr$:
\blr{5}, \blr{7}, and \blr{9}, respectively. 
Two sources are considered in the experiment 
shown in \cref{figure:wavefields-1mHz-fr-blr}, 
both located on the $z$-axis, at heights 
\num{0.7} and \num{0.99986}. 
The figure shows the solutions for a fixed 
radius of  $r=1$, slices on the 
meridional $xz$-plane, and the relative 
difference $\err_{xz}$, \cref{relative-difference:xz}.

By comparing the solutions shown in \cref{figure:wavefields-1mHz-fr-blr}, 
we can visually identify inaccuracies associated with the low-rank 
compression \blr{5} on the $xz$-plane slices. At this compression level, 
the solution loses its symmetric (with respect to the $z$-axis) structure, 
particularly for the source located near the surface. 
In contrast, with \blr{7}, no visible differences 
can be observed compared to the full-rank solution.
The maps of the relative difference confirm these observations: 
the low-rank factorization with threshold \blr{5} has significant 
differences, whereas the other thresholds yield 
solutions that match the full-rank ones. In these cases, the discrepancies are 
primarily localized  around the source. However, when the source is positioned 
near the surface, the differences extend farther, forming an inaccurate layer
near the surface, in particular for the low-rank \blr{7}.


We next evaluate the computational costs in \cref{table:blr-facto-mesh06Mb},
using the indicators of \cref{table:mumps-analytics}, 
in particular \pctNop, \pctNen~and~\pctNMP. 
Let us illustrate how storage of the factors is quantified with a simple example:
suppose that \fr~factorization needs to store 100 entries, 
while the BLR one needs 20: the storage is reduced to \pctNen$=20$\si{\percent}. 
Next, with mixed precision activated, say that 10 of these 
20 entries can be stored in half-precision instead of single, 
this gives $10 + 10\times0.5 = 15$ equivalent-entries
to store, and \pctNMP=$75$\si{\percent} ($15/20$). 
Combining BLR and mixed precision leads 
to \pctNen$\times$\pctNMP$=15\si{\percent}$ 
of entries compared to \fr~factorization. 

\begin{table}[ht!] \begin{center}
\caption{Computational cost for the matrix factorization with 
         \meshA~at frequency 1\si{\milli\Hz} 
         (solutions shown in \cref{figure:wavefields-1mHz-fr-blr})
         depending on the 
         BLR compression and the computational resources
         used. The matrix is of size $N=\num{132231488}$.
         The indicators are detailed in \cref{table:mumps-analytics}.}
\vspace*{-0.50em} \scriptsize
\label{table:blr-facto-mesh06Mb}
\renewcommand{\arraystretch}{1.20}
\sisetup{
exponent-product = \hspace*{-0.20em}\times\hspace*{-0.20em}, 
round-mode = places, 
round-precision = 2}

\begin{tabular}{|>{\arraybackslash}p{.25\linewidth}|
                 >{\arraybackslash}p{.12\linewidth}|
                 >{\arraybackslash}p{.12\linewidth}|
                 >{\arraybackslash}p{.12\linewidth}|
                 >{\arraybackslash}p{.12\linewidth}|
               }
\hline
  \multicolumn{5}{|c|}{Simulations with \meshA~at 1\si{\milli\Hz} using 
                       50 \adastra~nodes with 4MPI and 48threads per MPI.} \\ \hline
  & \texttt{full-rank} & \blr{5} & \blr{7} & \blr{9} \\ \hline
analysis time & \num{73.3981}\si{\second} 
              & \num{80.2753}\si{\second} 
              & \num{80.8232}\si{\second} 
              & \num{84.3185}\si{\second} \\ 
factorization memory & \num{12.432703}\si{\tera\byte}
                     & \num{6.6089360}\si{\tera\byte}
                     & \num{6.7043080}\si{\tera\byte}
                     & \num{6.8439090}\si{\tera\byte} \\ 
factorization time   & 1\si{\hour}03\si{\minute}
                     & 10\si{\minute}
                     & 12\si{\minute}
                     & 16\si{\minute} \\ 
 \texttt{RINFOG(3)}   & \num{106.4375439790e10} & - & - & - \\ 
\pctNop & 100\si{\percent}
                    & \num{ 2.9}\si{\percent}
                    & \num{ 6.1}\si{\percent}
                    & \num{10.6}\si{\percent} \\ 
\texttt{INFOG(29)}  & \num{105.9014083896e10} & - & - & - \\ 
\pctNen & 100\si{\percent}
                    & \num{13.1}\si{\percent}
                    & \num{21.6}\si{\percent}
                    & \num{30.3}\si{\percent} \\ 
\pctNMP & -
                    & \num{85.514582860821150}\si{\percent}
                    & \num{87.868109476582831}\si{\percent}
                    & \num{91.070151319162051}\si{\percent} \\ \hline
\end{tabular}
\vspace*{1em}

\begin{tabular}{|>{\arraybackslash}p{.25\linewidth}|
                 >{\arraybackslash}p{.12\linewidth}|
                 >{\arraybackslash}p{.12\linewidth}|
                 >{\arraybackslash}p{.12\linewidth}|
                 >{\arraybackslash}p{.12\linewidth}|
               }
  \hline
  \multicolumn{5}{|c|}{Simulations with \meshA~at 1\si{\milli\Hz} using 
                       20 \adastra~nodes with 4MPI and 48threads per MPI.} \\ \hline
  & \texttt{full-rank} & \blr{5} & \blr{7} & \blr{9} \\ \hline
analysis time & \texttt{x}
              & \num{75.8832}\si{\second} 
              & \num{76.2522}\si{\second} 
              & \num{77.2703}\si{\second} \\
factorization memory & \texttt{x}
                            & \num{4.377209}\si{\tera\byte}
                            & \num{5.239154}\si{\tera\byte}
                            & \num{5.822692}\si{\tera\byte} \\ 
factorization time   & \texttt{x}
                            & 10\si{\minute}
                            & 12\si{\minute}
                            & 16\si{\minute} \\
\pctNop & \texttt{x}
                    & \num{ 2.1}\si{\percent}
                    & \num{ 5.0}\si{\percent}
                    & \num{ 9.1}\si{\percent} \\
\pctNen & \texttt{x}
                    & \num{13.5}\si{\percent}
                    & \num{21.8}\si{\percent}
                    & \num{30.8}\si{\percent} \\ 
\pctNMP & \texttt{x}
                    & \num{85.759523724856322}\si{\percent}
                    & \num{87.704656681338022}\si{\percent}
                    & \num{91.756979899843017}\si{\percent} \\ \hline
\end{tabular}
\end{center}
\end{table}

From \cref{table:blr-facto-mesh06Mb}, 
we see that the BLR factorization leads to a drastic decrease in 
computational cost, both in terms of memory consumption and runtime. 
Due to the high-memory requirements, 
\fr~factorization need to use at least 50 nodes, while the BLR ones
can be executed on 20 nodes.
With 50 nodes, 
the BLR factorizations reduce the memory usage by 
approximately 55\si{\percent} and the factorization time 
by 80\si{\percent}. 
More precisely, using BLR factorization, the number of 
operations is in between 2 and 10\si{\percent} of the 
FR factorization,
and the number of entries, combining BLR and mixed-precision
storage (i.e., the product between \pctNop~and \pctNen), 
falls in between 11 and 28\si{\percent}, respectively
for \blr{5} and \blr{9}. 
Naturally, we see that increasing the BLR threshold reduces the 
cost, but in our case \blr{5} does not give accurate
solution (cf. \cref{figure:wavefields-1mHz-fr-blr}), and is
avoided in the following experiments. 
On the other hand, \blr{7} and \blr{9} already provide
significantly more frugal computations while maintaining 
a high-accuracy solution, as illustrated 
in \cref{figure:wavefields-1mHz-fr-blr}. 

Running the low-rank compression on 20 nodes leads to an additional 
reduction in memory usage, with comparable time costs, suggesting 
that allocating 50 nodes is excessive in this case, as our approach 
to parallelism is not strongly scalable in the sense that for a 
given problem size, increasing parallelism introduces 
overheads due to extra memory allocation for communication buffers.
While the differences between compression levels are minor when 
using a large number of nodes, they become clear with 20 nodes: 
the level \blr{7} further reduces memory consumption by 10\si{\percent} 
compared to \blr{9}.
We also notice the impact of the mixed precision with adaptive precision representation of the factors, 
that consistently contributes to about 10\si{\percent} of the 
reduction in the storage of the factors.

\medskip

The statistics for meshes \meshB~and~\meshC~are given 
in \cref{table:blr-facto-mesh08Mb}, for frequencies 
\num{1} and \num{1.5}\si{\milli\Hz}.
Here, the full-rank factorization was unable to run and we focus
on \blr{7} and \blr{9} which give accurate enough solutions.
Between the levels of compression, we have results that are 
consistent with \meshA~(\cref{table:blr-facto-mesh06Mb}) and
a more effective \blr{7} that gives between 
14 and \num{18}\si{\percent} memory reduction compared to \blr{9}.
We further note that the numbers of operations and entries are
only a fraction of an hypothetical FR factorization, and that 
the mixed-precision storage still contributes to an extra 
10\si{\percent} reduction.

\begin{table}[ht!] \begin{center}
\caption{Computational cost for the matrix factorization with 
         \meshB~and \meshC. The sizes of the resulting matrices
         with varying frequencies are given in \cref{table:matrix-analytics}.}
\vspace*{-0.50em} \scriptsize
\label{table:blr-facto-mesh08Mb}
\renewcommand{\arraystretch}{1.20}
\sisetup{
exponent-product = \hspace*{-0.20em}\times\hspace*{-0.20em}, 
round-mode = places, 
round-precision = 2}

\begin{tabular}{|>{\arraybackslash}p{.24\linewidth}|
                 >{\arraybackslash}p{.10\linewidth}|
                 >{\arraybackslash}p{.10\linewidth}||
                 >{\arraybackslash}p{.10\linewidth}|
                 >{\arraybackslash}p{.10\linewidth}|
               }
  \hline
  \multicolumn{5}{|c|}{Simulations with \meshB~using 
                       \adastra~nodes with 4MPI and 48threads per MPI.} \\ \hline
  & \multicolumn{2}{ c||}{ 1\si{\milli\Hz} with 30 nodes}
  & \multicolumn{2}{ c|}{1.5\si{\milli\Hz} with 50 nodes} \\ \hline
  & \blr{7} & \blr{9} & \blr{7} & \blr{9} \\ \hline
analysis time 
              & \num{103.1630}\si{\second}
              & \num{103.2348}\si{\second} 
              & \num{119.5584}\si{\second} 
              & \num{117.2157}\si{\second} \\
factorization memory 
                            & \num{7.405760}\si{\tera\byte}
                            & \num{8.608619}\si{\tera\byte} 
                            & \num{11.661180}\si{\tera\byte}
                            & \num{13.863476}\si{\tera\byte} \\ 
factorization time   
                            & 17\si{\minute}
                            & 22\si{\minute} 
                            & 27\si{\minute}  
                            & 41\si{\minute}  \\
\pctNop 
                    & \num{4.4}\si{\percent}
                    & \num{7.9}\si{\percent} 
                    & \num{7.1}\si{\percent}
                    & \num{12.7}\si{\percent} \\
\pctNen
                    & \num{18.9}\si{\percent}
                    & \num{26.8}\si{\percent}
                    & \num{21.2}\si{\percent}
                    & \num{30.7}\si{\percent} \\
\pctNMP           
                    & \num{87.407025742857144}\si{\percent}
                    & \num{91.172129218177488}\si{\percent}
                    & \num{85.945173094546959}\si{\percent} 
                    & \num{90.229684399050910}\si{\percent} \\ \hline
\end{tabular}

\vspace*{1em}

\begin{tabular}{|>{\arraybackslash}p{.30\linewidth}|
                 >{\arraybackslash}p{.20\linewidth}|
                 >{\arraybackslash}p{.20\linewidth}|
               }
  \hline
  \multicolumn{3}{|c|}{Simulations with \meshC~using 80 \adastra~nodes with 4MPI and 48threads per MPI.} \\ \hline
  & \blr{7} -- 1\si{\milli\Hz}
  & \blr{7} -- 1.5\si{\milli\Hz} \\ \hline
analysis time & \num{161.4035}\si{\second}
              & \num{169.6870}\si{\second} \\
factorization memory & \num{16.039254}\si{\tera\byte} &  \num{19.624699}\si{\tera\byte} \\ 
factorization time   & 27\si{\minute}  &  42\si{\minute} \\
\pctNop & \num{4.8}\si{\percent}       & \num{6.5}\si{\percent} \\
\pctNen 
                    & \num{17}\si{\percent} & \num{19}\si{\percent} \\
\pctNMP & \num{87.209540404749205}\si{\percent}
        & \num{85.341554761069347}\si{\percent} \\ \hline
\end{tabular}

\end{center}
\end{table}

Comparing the different meshes, we see that the computational
cost increases with the more refined mesh, in terms of memory
consumption, time, and resources required to run. For instance
at 1\si{\milli\Hz} with \blr{7}, the factorization memory is 
of \num{5.24}, \num{7.41} and \num{16.04}\si{\tera\byte} for
\meshA, \meshB~and~\meshC, respectively, i.e., an 
increase by a factor 3 for the finer mesh compared to the coarser one. 
To maintain accuracy, the coarse meshes compensate with 
higher-order polynomials, which still gives a much smaller linear 
system than with refinement (see \cref{table:matrix-analytics}), 
thus they retain a relatively low computational cost.
In \cref{table:blr-facto-2mHz}, we further compare the matrix 
factorization at 2\si{\milli\Hz} for the three meshes using \blr{7}. 
Here, the differences between the meshes are reduced. For example the factorization memory at 1\si{\milli\Hz} 
increases by a factor three between the coarser mesh \meshA ~and the finer mesh \meshC ~whereas at 2\si{\milli\Hz} it is slight less than a factor of two increase.
These observations on computational cost must now be considered 
in relation with the accuracy of the solutions, which is studied
in the following \cref{section:validation-with-axisym}.

\begin{table}[ht!] \begin{center}
\caption{Computational cost at 2\si{\milli\Hz} using \blr{7}.}
\vspace*{-0.50em} \scriptsize
\label{table:blr-facto-2mHz}
\renewcommand{\arraystretch}{1.20}
\sisetup{
exponent-product = \hspace*{-0.20em}\times\hspace*{-0.20em}, 
round-mode = places, 
round-precision = 2}

\begin{tabular}{|>{\arraybackslash}p{.24\linewidth}|
                 >{\arraybackslash}p{.15\linewidth}|
                 >{\arraybackslash}p{.15\linewidth}|
                 >{\arraybackslash}p{.15\linewidth}|
               }
\hline
 &   \meshA~using \newline 75 \adastra~nodes 
 &   \meshB~using \newline 75 \adastra~nodes 
 &   \meshC~using \newline 100 \adastra~nodes  \\ \hline
matrix size  & \num[round-precision=0]{202431918} & \num[round-precision=0]{257576619} & \num[round-precision=0]{309859823} \\ \hline
graph  size  & \num[round-precision=0]{12591470}  & \num[round-precision=0]{17097752}  &  \num[round-precision=0]{23112716} \\
edge   size  & \num[round-precision=0]{75182304}  & \num[round-precision=0]{102219996} & \num[round-precision=0]{138309780} \\
analysis time& \num{102.9651}\si{\second} & \num{138.9138}\si{\second} & \num{171.5996}\si{\second} \\ \hline
factorization memory         & \num{15.226126}\si{\tera\byte} & 
\num{19.893622}\si{\tera\byte} & \num{28.609103}\si{\tera\byte} \\ 
factorization time           & 38\si{\minute}               & 53\si{\minute}               & 54\si{\minute}               \\
\pctNop & \num{11.2}\si{\percent}      & \num{9.4}\si{\percent}       & \num{7.5}\si{\percent}       \\
\pctNen & \num{26.7}\si{\percent}      & \num{23.2}\si{\percent}      & \num{18.6}\si{\percent}      \\ 
\pctNMP & \num{84.61930984656735}\si{\percent}
                              & \num{84.99699209842187}\si{\percent}
                              & \num{93.67596820009422}\si{\percent} \\ \hline
\end{tabular}
\end{center}
\end{table}

\section{Validation with 2.5D solver}
\label{section:validation-with-axisym}

With helioseismology as the target application, the validation of the 
solver has to be carried out in relevant settings concerning particularly background 
parameters and source and receiver positions.
Due to the heterogeneous nature of the Sun's interior (\cref{figure:modelS}), 
physics-pertinent analytical solutions are not available. 
As references, we employ solutions
constructed in \cite{Pham2025axisymreport} via an azithmuthal expansion (2.5D)
for axially symmetric solar backgrounds.
With rotational symmetry around $\mathbf{e}_z$-axis and a scalar Dirac 
source along $\mathbf{e}_z$, the restriction to any meridional 
half-plane of a 3D solution coincides with the 2.5D solution 
at azithmuthal order $0$. 
The 2.5D implementation has the additional virtue of being able to use
a more refined mesh compared to 3D. This is particularly important 
in a neighborhood of Dirac sources, and in the near-surface layer 
where heights of excitation (i.e., source position) 
and observation (i.e., receiver position) are taken in helioseismology.

In the following, we consider the standard solar model S,
\cref{subsection:strategy-mesh}. 
We first compare solutions on the meridional half-plane
in \cref{mercom::subsec}, followed by comparisons at 
specific heights in \cref{fixheicomp::subsec}. 
As representatives of meridional half-plane, we consider those in
the $xz$-plane, at either azithmuthal angle $\phi = 0$ or $\pi$.
In these experiments, the attenuation is fixed to 
$\gamma_{\mathrm{att}}/(2\pi)=$10\si{\micro\Hz}, and 
in accordance with the findings in \cref{section:mumps}, 
$\blr{7}$ is employed within \mumps. 
We investigate frequencies 1\si{\milli\Hz} and 2\si{\milli\Hz}, 
and 3D meshes, \meshA,~\meshB, and \meshC~introduced in \cref{section:implementation}.
Regarding source and receiver configurations, we
consider three Dirac source positions on the positive 
$z$-axis 
and situated at heights: $\zsrc=0.7$, $\zsrc=0.99$ and $\zsrc=0.99986$. 
The first position is situated in the transition layer 
between the radiative and the convective zone, while the 
last two are in the near-surface layer, 
respectively at \num{6960} and  \num{97}\si{\kilo\meter} below surface.

The frequencies \num{1}\si{\milli\Hz} and \num{2}\si{\milli\Hz} 
are below the acoustic cut-off frequency ($\sim \num{5.3}\si{\milli\Hz}$)
and the computational domain consists of an interior elliptic 
region $(\omega > N)$ and a hyperbolic region  $(\omega < N)$ 
situated in the near-surface layer, \cite{Pham2025axisymreport}. 
As shown in \cite[Section 6]{Pham2025axisymreport}, for such frequencies, 
even with attenuation, achieving accuracy for solutions 
in the solar atmosphere is particularly challenging when a singular source is 
located in the near-surface layer.

\subsection{Comparison of the solutions on a meridional half-plane}
\label{mercom::subsec}

\begin{figure}[ht!]\centering

\subfloat[1\si{\milli\Hz} solution at $r=\num{1}$ 
          for a source at $\zsrc=\num{0.99986}$.]{\raisebox{1.5em}{
  \includegraphics[height=3.2cm]{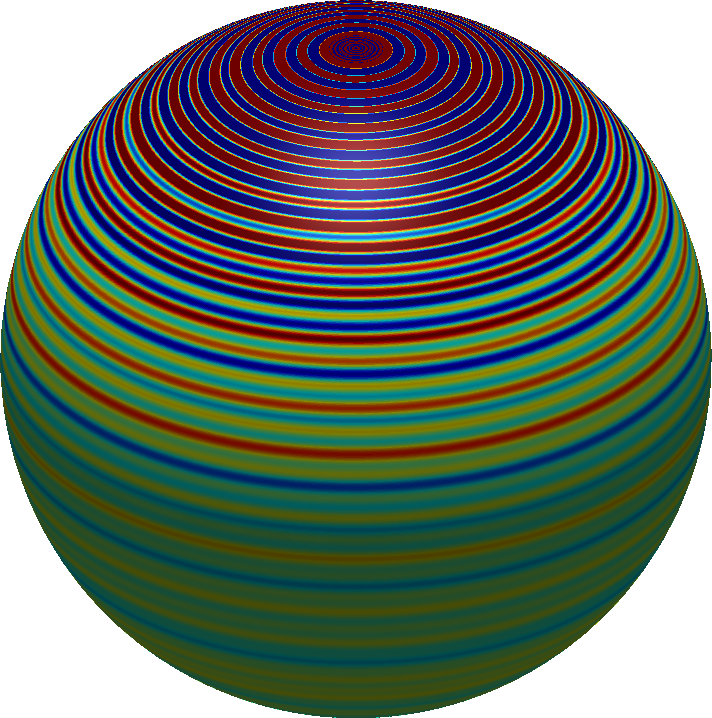}
          
  \pgfkeys{/pgf/fpu=true} 
  \pgfmathsetmacro{\cmin}{-5e-8} 
  \pgfmathsetmacro{\cmax}{5e-8} 
  \pgfkeys{/pgf/fpu=false}
  \setlength{\cbarwidth} {0.50em}
  \setlength{\cbarheight}{2.50cm}
  \hspace*{-4em}\raisebox{1.2em}{\input{figures/cbar}}
 
          \label{afigure:validation-mesh11M_1mHz-2mHz}}
\hspace*{-4em}{\raisebox{1.0em}{\begin{tikzpicture}
  \node[] at (0,0) {\scriptsize 1\si{\milli\Hz}} ; 
\end{tikzpicture}}}
} \hspace*{1em}
\subfloat[Solutions on the $xz$-meridional 
          half-plane at 1\si{\milli\Hz}
          for two source heights.]{
  \includegraphics[scale=1]{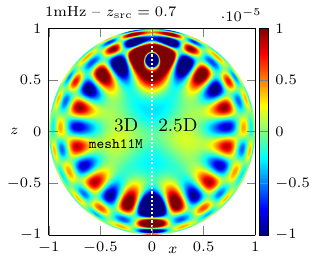}
  \includegraphics[scale=1]{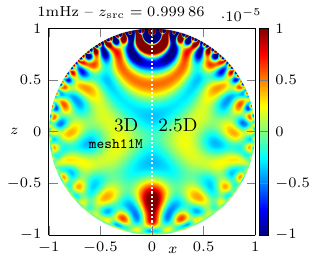}  
        \label{bfigure:validation-mesh11M_1mHz-2mHz}  
  }
\vspace*{-1em}

\subfloat[2\si{\milli\Hz} solution at $r=\num{1}$ 
          for a source at $\zsrc=\num{0.99986}$.]{\raisebox{1.5em}{
  \includegraphics[height=3.2cm]{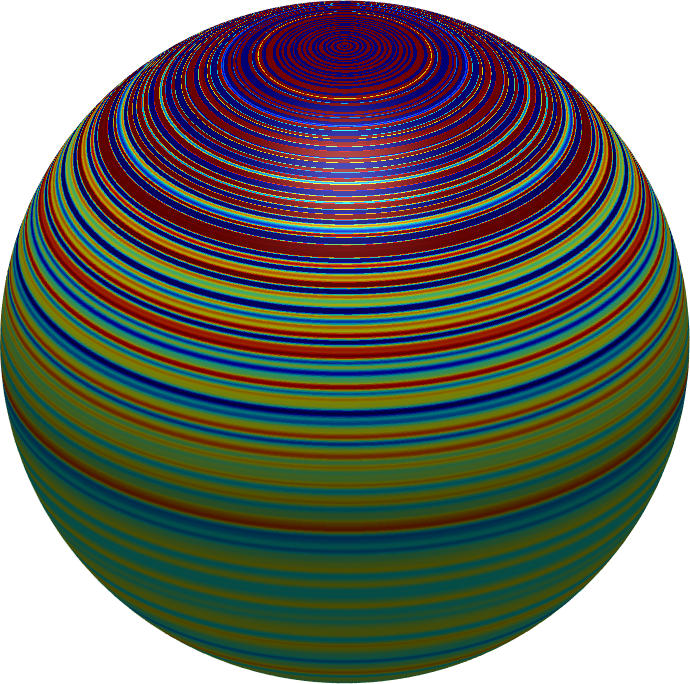}
          
  \pgfkeys{/pgf/fpu=true} 
  \pgfmathsetmacro{\cmin}{-4e-6} 
  \pgfmathsetmacro{\cmax}{4e-6} 
  \pgfkeys{/pgf/fpu=false}
  \setlength{\cbarwidth} {0.50em}
  \setlength{\cbarheight}{2.50cm}
  \hspace*{-4em}\raisebox{1.2em}{\input{figures/cbar}}

          \label{cfigure:validation-mesh11M_1mHz-2mHz}}
\hspace*{-4em}{\raisebox{1.0em}{\begin{tikzpicture}
  \node[] at (0,0) {\scriptsize 2\si{\milli\Hz}} ; 
\end{tikzpicture}}}} \hspace*{1em}
\subfloat[Solutions on the $xz$-meridional 
          half-plane at 2\si{\milli\Hz}
          for two source heights.]{
  \includegraphics[scale=1]{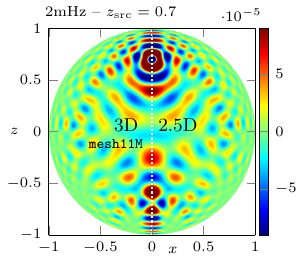}
  \includegraphics[scale=1]{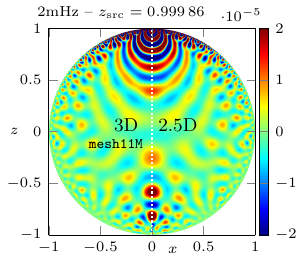}  
  \label{dfigure:validation-mesh11M_1mHz-2mHz}}

\caption{Comparison of the 3D and 2.5D simulations at frequency 1\si{\milli\Hz} 
         (\textbf{top}) and 2\si{\milli\Hz} (\textbf{bottom}) with the solar 
         background models.
         The $xz$-meridional half-plane solutions in (\textbf{b}, \textbf{d}) 
         show the 3D simulations on the left, and 2.5D ones on the right.
         The 3D simulations use \meshC~and \blr{7}. 
         }
\label{figure:validation-mesh11M_1mHz-2mHz}
\end{figure}

We first compare solutions on the $xz$-meridional half-plane.
Solutions at frequencies 1\si{\milli\Hz} and 2\si{\milli\Hz} are shown 
in \cref{figure:validation-mesh11M_1mHz-2mHz}, with \meshC~employed for 
the 3D resolutions. 
Two Dirac source positions along the $z$-axis are considered 
in \cref{figure:validation-mesh11M_1mHz-2mHz}: 
at $\zsrc = 0.7$ (top row), and $\zsrc = 0.99986$ (bottom row).
The left panels, \cref{afigure:validation-mesh11M_1mHz-2mHz,cfigure:validation-mesh11M_1mHz-2mHz}, 
show the solution at a fixed radial height $r = 1$, 
and feature spherical wave pattern.
Here, no significant differences can be observed between the 
3D and 2.5D solutions, which exhibit identical amplitudes and 
oscillation pattern in 
\cref{bfigure:validation-mesh11M_1mHz-2mHz,dfigure:validation-mesh11M_1mHz-2mHz}.

\begin{figure}[ht!] \centering
  \subfloat[3D with \meshA.]
           {\includegraphics[scale=1]{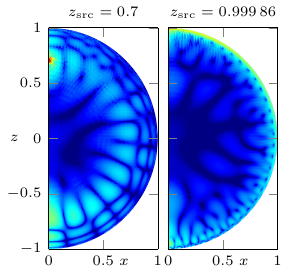}} \hfill
  \subfloat[3D with \meshB.]{
              \includegraphics[scale=1]{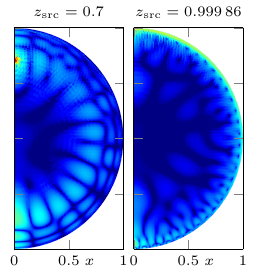}} \hfill
  \subfloat[3D with \meshC.]{
              \includegraphics[scale=1]{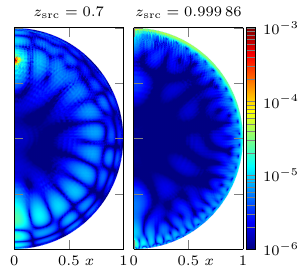}}

\caption{Comparison of the relative difference $\err_{xz}$ \cref{relative-difference:xz}
         at frequency 1\si{\milli\Hz} for various the source positions and 3D meshes. 
         The reference 2.5D solution and the 3D one using \meshC~are shown in
         \cref{figure:validation-mesh11M_1mHz-2mHz}.}
\label{figure:validation-relative-diff-map}
\end{figure}

To better detect differences in \cref{figure:validation-relative-diff-map}, 
we map the relative difference $\err_{xz}$ defined in \cref{relative-difference:xz},
between 3D solutions obtained on the three meshes and the 2.5D reference solution.
The comparison is for 1\si{\milli\Hz} and the two Dirac source positions. 
We have the following observations.

\begin{itemize}[leftmargin=*]\setlength{\itemsep}{-1pt}
  \item All three 3D meshes yield comparable relative difference.
        This highlights the effectiveness of the numerical strategy introduced 
        in \cref{section:implementation}: coarser meshes are 
        compensated by higher-order polynomial functions, 
        with a cell-wise accurate representation of background 
        model as described in \cref{remark:model-representation}. 
        Comparatively, the biggest difference is observed with \meshA, 
        while \meshB~and \meshC~yield similar relative differences.
  \item Most prominent differences are observed in an annulus layer 
        at the source height, reaching a maximum around the source 
        position $\zsrc$. The differences decrease with distance 
        to the source, with the exception of an increase at $-\zsrc$.
\end{itemize}
The above experiment demonstrates the well-known fact that
maintaining accuracy near a Dirac source position is challenging 
and would require an extremely fine mesh at this position. 
Futhermore, in our cases, the use of \blr{7} in the 3D simulations 
leads to an additional drop in accuracy in this region, 
as already observed in \cref{figure:wavefields-1mHz-fr-blr}.

\subsection{Comparison of solutions at fixed heights}
\label{fixheicomp::subsec}

In computing synthetic helioseismic observables,
single-height observation assumption is typically used, which means that
only the trace of the numerical solutions at a fixed height $r$ is employed.
In view of this, we compare values of solutions at 
 $r=1$ and $r=0.7$. The height 
$r=1$, at the surface, constitutes the observation height employed in \cref{section:numerics-3d}, 
while the height $r=0.7$, situated below the convective zone and thus away from the near-surface layer,
offers a comparative contrast.
We introduce the point-wise relative error $\err_r$
defined for a fixed height as,
\begin{equation}\label{eq:relative-difference:r}
  \err_{r}(\theta) \,=\, 
  \dfrac{\vert w_{\texttt{2.5D}}(\bx) - w_{\texttt{3D}}(\bx) \vert}{\Vert w_{\texttt{2.5D}} \Vert} 
  \,, \hspace*{2em} 
  \text{where $\bx = \left(r \cos(\theta), 0, r \sin(\theta)\right)$, 
  \quad $\theta \in (-\pi/2,\pi/2)$.}
\end{equation}
Here, $w_{\texttt{2.5D}}$ is the solution computed with the axisymmetric
solver and $w_{\texttt{3D}}$ the 3D solution.
The Dirac source is placed along the $z$-axis, thus corresponding to $\theta = \pi/2$.
We consider an interior Dirac source at $\zsrc=0.7$,
and for the near-surface layers, $\zsrc=0.99$ and $\zsrc=0.99986$.

\subsubsection{Interior source} 

\begin{figure}[ht!] \centering
\subfloat
{\includegraphics[scale=1]{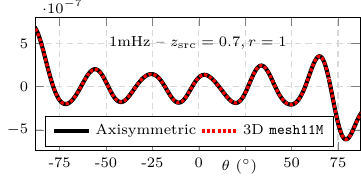} \hspace*{1em}
 \includegraphics[scale=1]{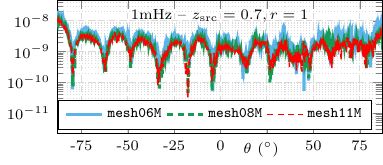}}
\vspace*{-0.75em}

\subfloat
{\includegraphics[scale=1]{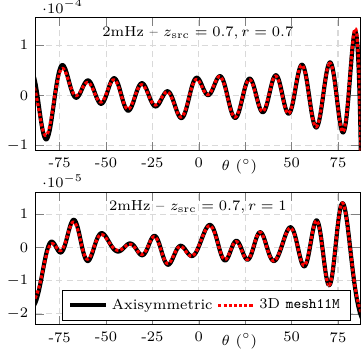} \hspace*{1em}
 \includegraphics[scale=1]{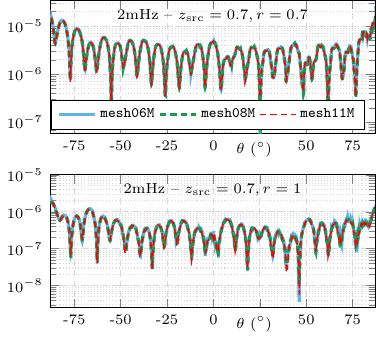}}
 
\caption{Comparisons of solutions at fixed height $r$ for a source at $\zsrc=\num{0.7}$ for
         frequencies 1\si{\milli\Hz} and 2\si{\milli\Hz}. The wavefields are shown on the {left}
         and their relative difference $\err_r$ \cref{eq:relative-difference:r} with the 
         2.5D solutions on the {right}.}
\label{figure:validation-radial-interior-src}
\end{figure}

We first compare the solutions corresponding to a Dirac source 
located at $\zsrc= \num{0.7}$ (recall that wavefields were shown 
in the middle column of \cref{figure:validation-mesh11M_1mHz-2mHz}).
The left column of \cref{figure:validation-radial-interior-src} 
shows cuts of the solutions at fixed height $r=1$ and $r=0.7$, 
for frequencies 1\si{\milli\Hz} and 2\si{\milli\Hz}.
The right column compares solutions in terms of the relative 
differences $\err_r$ \cref{eq:relative-difference:r} among the 
meshes \meshA,~\meshB, and \meshC~introduced in \cref{section:implementation}.

In this case, we observe an agreement 
between the 2.5D and 3D solutions, with no visible difference between
in the wavefields (left of \cref{figure:validation-radial-interior-src}), 
and the relative difference remains below $\num{e-5}$, 
irrespective of the receiver height.
All three 3D meshes yield the same level of relative difference,
in accordance with \cref{figure:validation-relative-diff-map}.

\subsubsection{Near-surface sources} 

\begin{figure}[ht!] \centering
\subfloat
{\includegraphics[scale=1]{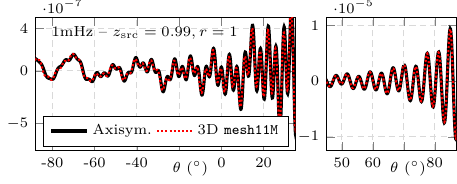}  \hspace*{.2em}
 \includegraphics[scale=1]{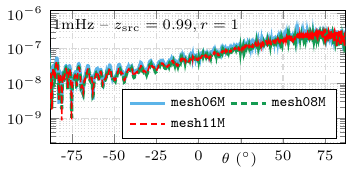}}
\vspace*{-1em}

\subfloat
{\includegraphics[scale=1]{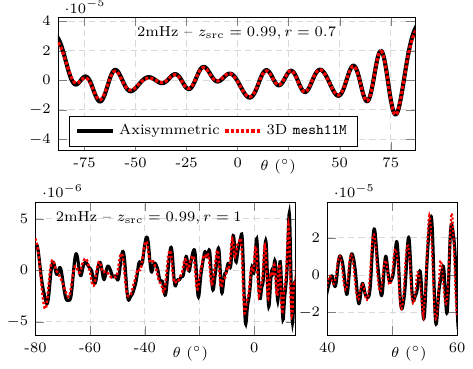} \hspace*{.2em}
 \includegraphics[scale=1]{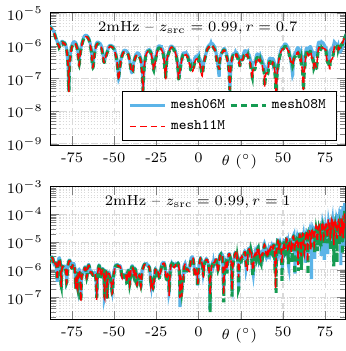}}
 
\caption{Comparisons of solutions at fixed height $r$ for a source at $\zsrc=\num{0.99}$ for
         frequencies 1\si{\milli\Hz} and 2\si{\milli\Hz}. The wavefields are shown on the {left}
         and their relative difference $\err_r$ \cref{eq:relative-difference:r} with the 
         2.5D solutions on the {right}. At height $r=1$, two angular regions are shown
         on different scales for better visualization.}
\label{figure:validation-radial-exterior-src0.99}
\end{figure}

\begin{figure}[ht!] \centering
\subfloat
{\includegraphics[scale=1]{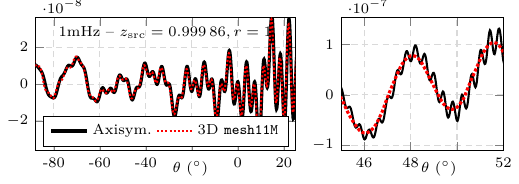}  \hspace*{.2em}
 \includegraphics[scale=1]{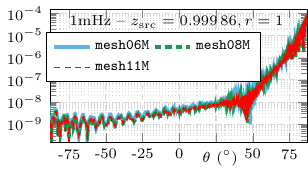}}
\vspace*{-1em}

\subfloat
{\includegraphics[scale=1]{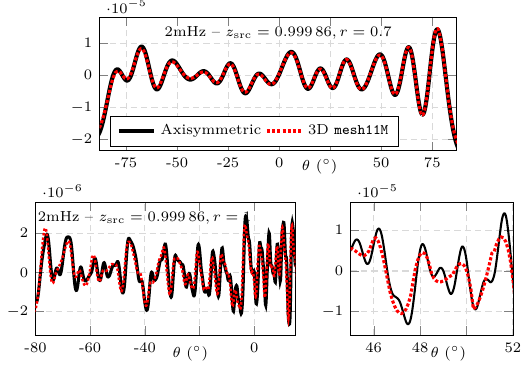} \hspace*{.2em}
 \includegraphics[scale=1]{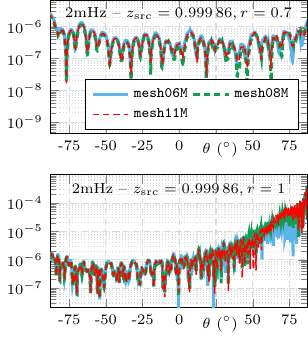}}
 
\caption{Comparisons of solutions at fixed height $r$ for a source at $\zsrc=\num{0.99986}$ for
         frequencies 1\si{\milli\Hz} and 2\si{\milli\Hz}. The wavefields are shown on the {left}
         and their relative difference $\err_r$ \cref{eq:relative-difference:r} with the 
         2.5D solutions on the {right}. At height $r=1$, two angular regions are shown
         on different scales for better visualization.}
\label{figure:validation-radial-exterior-src0.99986}
\end{figure}

We now consider the case of a source located in the near-surface layer, 
at heights $\zsrc= \num{0.99}$ and $\zsrc= \num{0.99986}$. The 
corresponding wavefields and relative differences are shown 
in \cref{figure:validation-radial-exterior-src0.99,figure:validation-radial-exterior-src0.99986}, respectively.
Due to the rapid oscillations in the surface layers, the plots are separated
into two angular regions, near the source and away from it, 
featuring different scales to enhance visualization.
Recall that the Dirac source is always located at $\theta=90\si{\degree}$.
As discussed in \cite{Pham2025axisymreport}, this configuration 
is particularly challenging due to the change in the nature of 
the PDE in the near-surface layer at vanishing attenuation, 
which occurs for frequencies below the acoustic cut-off.
Nevertheless, the setup with both sources 
and receivers situated in the near-surface layer
is fundamental, as it is employed to generate 
synthetic observables in helioseismic applications.

For the source at $\zsrc=0.99$ (\cref{figure:validation-radial-exterior-src0.99}), 
the solution at height $r=0.7$ is accurate, showing only small 
differences compared to the 2.5D solution. 
However, at the surface height $r=1$, the 3D solution shows 
some discrepancies at 2\si{\milli\Hz}, mostly in amplitude. 
For both 1\si{\milli\Hz} and 2\si{\milli\Hz} frequencies, 
the relative difference increases as the receiver position approaches the source, 
particularly within the $45\si{\degree}$ of it. 
This indicates that the discrepancy is localized near the source and does 
not propagate into the interior. This behaviour 
is consistent with the observations in \cite{Pham2025axisymreport}, 
which reported on the confinement of singularity propagation 
in the near-surface layer, with strength concentrated in the 
neighborhood of the Dirac source.
In 3D, the difficulty to refine meshes near the singularity 
of the Dirac, and the use of BLR factorization further reduce the accuracy.

We observe a similar behaviour for the source at $\zsrc= 0.99986$ 
in \cref{figure:validation-radial-exterior-src0.99986}. The solution 
remains accurate in the interior, while discrepancies appear at the surface, 
especially at positions close to the source. 
Interestingly, at \num{1}\si{\milli\Hz}, the 3D solution near the source 
capture the main oscillations but not the finer-scale ones.
For both frequencies, within the $45\si{\degree}$ 
of the source, the relative differences increase by several orders of magnitude. 
This suggests that the 3D simulations can benefit from further local
refinement around the source position for higher accuracy.
For this reason, in \cref{section:numerics-3d}, to avoid complications 
due to Dirac singularities,we will consider a Gaussian source instead. 
\section{Solar wave propagation with 3D wave speed perturbations}
\label{section:numerics-3d}

In this section, we consider a non-spherical solar background 
with a 3D perturbation in wave speed such that,
\begin{equation}\label{eq:c3d-default}
  c_\delta(r,\theta,\phi) \,=\, c_0(r) + \delta_c(r,\theta,\phi) \,,
\end{equation}
where $c_0$ is the background wave speed of 
standard model S \cite{christensen1996current}
(see \cref{figure:modelS}), and $\delta_c$ is 
the 3D perturbation.
We design two configurations for the perturbations: one based on a map of active regions
in \cref{subsection:num-active-region}, and the other derived from a snapshot of a convection simulation,
\cref{subsection:num-convection}.
The effect of the 3D perturbations are emphasized by comparing
the difference with a solution computed with the radially symmetric
background $c_0$. This comparison is made at fixed height, and in the 
spirit of \cref{eq:relative-difference:r} we introduce the 
relative difference $\err_\mathfrak{r}$ such that,
\begin{equation}\label{eq:relative-difference:map}
  \err_\mathfrak{r}(\theta,\phi) \,=\, 
  \dfrac{\vert w_{c_0}(\bx) - w_{c_\delta}(\bx) \vert}{\Vert w_{c_0}(\bx) \Vert} 
  \,, \hspace*{2em} 
  \text{for $\bx=(r,\theta,\phi)$ with $r=\mathfrak{r}, \, \forall \theta\, \,\, \forall \phi$}.
\end{equation}
Here, $w_{c_0}$ is the solution computed with 
the spherically symmetric background of model S, 
and $w_{c_\delta}$ is the solution computed with 
the 3D perturbations \cref{eq:c3d-default}.

To avoid any numerical instabilities for a Dirac, 
we instead use a Gaussian source function for $f$ 
in \cref{eq:main-generic} such that,
\begin{equation}
  f(r,\theta,\phi) \, = \, \exp\left( -\,\dfrac{(r-r_s)^2}{\sigma_r^2} 
                                    \,-\,\dfrac{(\theta-\theta_s)^2}{\sigma_\theta^2} 
                                    \,-\,\dfrac{(\phi-\phi_s)^2}{\sigma_\phi^2} 
                               \right) \,,
\end{equation}
where $(r_s,\theta_s,\phi_s)$ is the center of the source. 
Motivated by the solar stratification, it is taken to be 
narrow in $r$ with variance $\sigma_r\,=\,\num{e-4}$, while 
$\sigma_\theta=\sigma_\phi=0.2$\si{\radian}.
The experiments are carried out using \meshB~with \blr{9} 
and the level attenuation is $\gamma_{\mathrm{att}}/(2\pi)=$10\si{\micro\Hz}.

\subsection{Wave propagation with solar active regions}
\label{subsection:num-active-region}

\subsubsection{3D wave speed perturbation from active regions}

\begin{figure}[ht!]\centering
\includegraphics[scale=1.1]{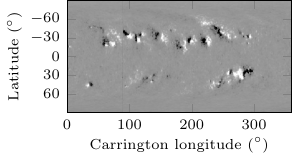} \hspace*{1em}
\raisebox{0.50em}{
\includegraphics[scale=1.1]{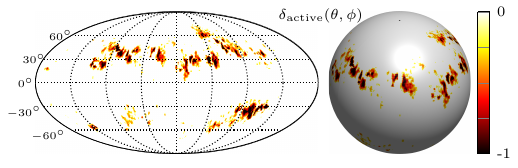}}
\caption{Synoptic magnetogram of the Carrington rotation 2108 
         from March 2011 (\textbf{left}).
         The magnetic activity is converted into a perturbation
         $\delta_{\mathrm{active}} \in (-1,0)$ 
         visualized as Mollweide map (\textbf{center}) 
         and on the sphere (\textbf{right}) at the height where
         it is maximal, $r_{\mathrm{active}}=0.995$. 
         It is then used to define the 3D wave speed in \cref{eq:wavespeed-active}.
         }
\label{fig:active-region:model}
\end{figure}

We create the 3D perturbation from the magnetogram synoptic map 
from GONG\footnote{\href{https://magmap.nso.edu/}{https://magmap.nso.edu/}}. 
An example from the Carrington rotation 2108 corresponding
to March 2011 is
shown in \cref{fig:active-region:model}. As the active regions are cooler than the quiet Sun, it results in a decrease of the wave speed in these
areas, hence $\delta_{\mathrm{active}} \in (-1,0)$,
cf. \cref{fig:active-region:model}. The 3D wave speed $c_{\mathrm{active}}$ is, 
\begin{equation}\label{eq:wavespeed-active}
  c_{\mathrm{active}}(r,\theta,\phi) \,=\, c_0(r) \,\,\left( \, 1 + \, 
     \alpha \, g(r) \, \delta_{\mathrm{active}}(\theta,\phi) \,\right) \,,
\end{equation}
where $g(r)$ is a Gaussian function for which we select the mean at height 
$r_{\mathrm{active}}=\num{0.995}$ and variance \num{8.5e-4}, i.e., it is 
relatively narrow in $r$. The parameter $\alpha\in(0,1)$ serves to adjust the 
maximal scale of the perturbation compared to the background.

\subsubsection{Effect of the perturbations on the wavefields}

We compare the wavefield solutions computed using the radially 
symmetric background $c_0$ and the perturbed 3D background $c_{\mathrm{active}}$. 
In \cref{fig:active-region:wavefields-20pct_srcpos}, we show the solutions 
at height $r = 1$ for a source located at 
$(r_s, \theta_s, \phi_s) = (0.99, \num{60}\si{\degree}, \num{135}\si{\degree})$, 
at frequency 1\si{\milli\Hz}. Here we use a scaling $\alpha=20\si{\percent}$ for
the 3D perturbation \cref{eq:wavespeed-active}.
The differences between the two solutions are further displayed on a logarithmic 
scale to enhance visibility.

\begin{figure}[ht!]\centering
\includegraphics[scale=1]{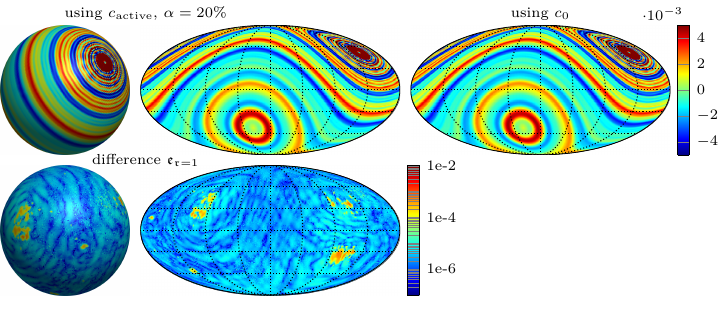}

\includegraphics[scale=1]{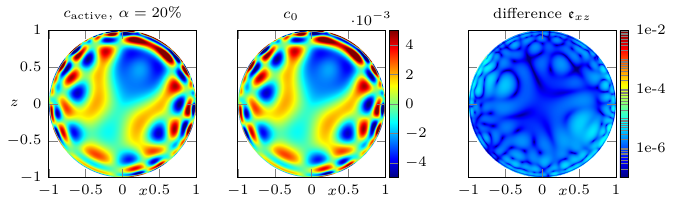}

\caption{Solutions for a source positioned at 
         $(r_s, \theta_s, \phi_s) = (0.99, \num{60}\si{\degree}, \num{135}\si{\degree})$
         for 1\si{\milli\Hz} frequency. 
         On top we show the solutions at fixed height $r=1$
         and on the bottom the solution on the meridional plane.
         The simulation with the 3D wave speed perturbation 
         \cref{eq:wavespeed-active} uses $\alpha=20\si{\percent}$.}
\label{fig:active-region:wavefields-20pct_1mHz}
\end{figure}

We observe that the solutions computed with the radial background $c_0$ 
and with the perturbed background $c_{\mathrm{active}}$ with $\alpha = 20\si{\percent}$
cannot be distinguished visually, both at the surface ($r = 1$) 
and in the interior slice.  
However, the maps of the relative difference reveal the pattern of the active regions shown 
in \cref{fig:active-region:model}, especially those with the highest amplitudes, which 
match the position of the largest differences. 
This is further confirmed in \cref{fig:active-region:wavefields-20pct_srcpos}, 
where we vary the source positions with height. In each case, 
the position of the strongest active regions perturbations 
correspond to the areas with the most pronounced differences 
between the solutions.
That is, the map of the difference is barely affected by the position 
of the source.
Nonetheless, for the highest source position ($s_r = \num{0.99986}$), we observe an
increase in difference near the source position. This behaviour is consistent with the 
observations drawn in \cref{section:validation-with-axisym}, with difficulties
handling sources too close to the surface, and because of inaccuracies induced 
by the use of block low-rank factorization.

\begin{figure}[ht!]\centering
\includegraphics[scale=1]{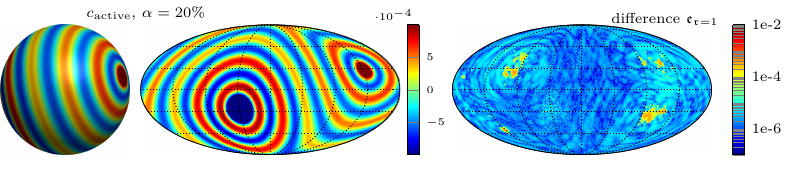}
\vspace*{-0.80em}

\includegraphics[scale=1]{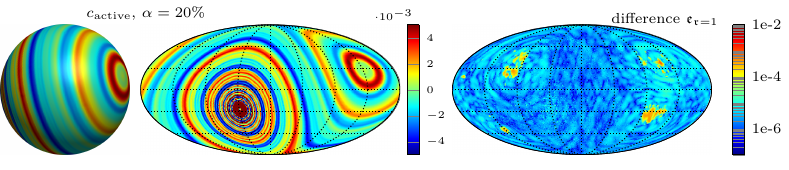}
\vspace*{-0.80em}

\includegraphics[scale=1]{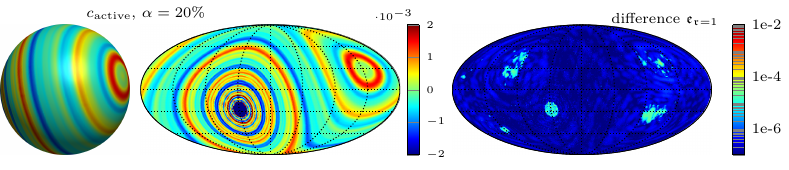}

\caption{Simulations at 1\si{\milli\Hz} for a source positioned at 
         $(\theta_s, \phi_s) = (\num{-22.5}\si{\degree}, \num{-45}\si{\degree})$
         and heights $r_s=$\num{0.7} (\textbf{top}), \num{0.99} (\textbf{middle})
         and \num{0.99986} (\textbf{bottom}). 
         The solutions and relative difference $\err_R$ are 
         pictured at height $r=1$.
         The simulation with the 3D wave speed perturbation 
         \cref{eq:wavespeed-active} uses $\alpha=20\si{\percent}$.}
\label{fig:active-region:wavefields-20pct_srcpos}
\end{figure}

\begin{figure}[ht!]\centering
\includegraphics[scale=1]{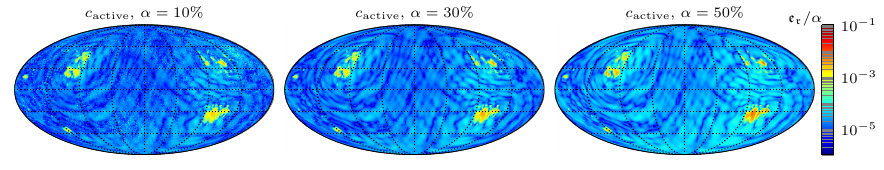}

\includegraphics[scale=1]{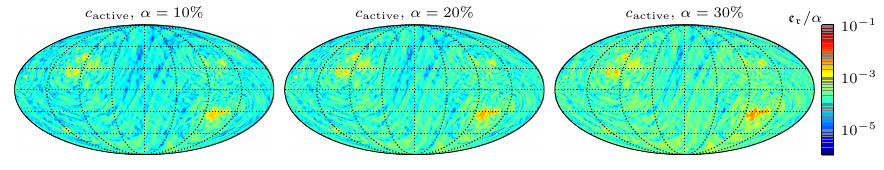}
\caption{Relative difference $\err_R/\alpha$ 
         at height $r=1$ at \num{1}\si{\milli\Hz} (\textbf{top})
         and \num{2}\si{\milli\Hz} (\textbf{bottom}) for different
         level of perturbation $\alpha$ \cref{eq:wavespeed-active}. 
         The source is positioned in 
         $(r_s, \theta_s, \phi_s) = (0.99, \num{60}\si{\degree}, \num{135}\si{\degree})$.}
\label{fig:active-region:wavefields-XXpct_1mHz}
\end{figure}

\begin{figure}[ht!]\centering
\includegraphics[scale=1]{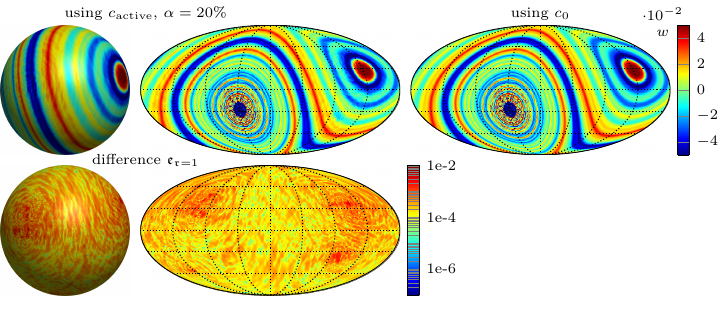}
\caption{Solutions at 2\si{\milli\Hz} for a source positioned at 
         $(r_s, \theta_s, \phi_s) = (0.99, \num{-22.5}\si{\degree}, \num{-45}\si{\degree})$,
         showing the height $r=1$.
         The simulation with the 3D wave speed perturbation 
         \cref{eq:wavespeed-active} uses $\alpha=20\si{\percent}$.}
\label{fig:active-region:wavefields-20pct_high-freq}
\end{figure}

In \cref{fig:active-region:wavefields-XXpct_1mHz}, we vary the perturbation 
amplitude $\alpha$ from 10 to 50\si{\percent} and present the corresponding 
difference maps scaled $\err_{\mathfrak{r}=1}/\alpha$, 
for frequencies \num{1} and \num{1.5}\si{\milli\Hz}. 
These experiments emphasize that the relationship between the wave speed and wavefields perturbations remain mostly linear, and thus the Born approximation is valid. Deviations start to be visible when the perturbation level reaches 50\%.
We further consider frequency of 2\si{\milli\Hz} and show 
the corresponding solutions with $\alpha = 20\si{\percent}$ 
in \cref{fig:active-region:wavefields-20pct_high-freq}.
Here, visual discrepancies can be seen in the wavefields 
with and without 3D perturbations. 
In the maps of the relative difference, the patterns of the active 
regions are no longer clearly visible, even though the highest 
differences still occur at their locations. 
Namely, the background difference is of the same magnitude
as the difference at the location of the active region.

\subsection{Wave propagation through convection}
\label{subsection:num-convection}

\begin{figure}[ht!]\centering
           {\includegraphics[scale=1]{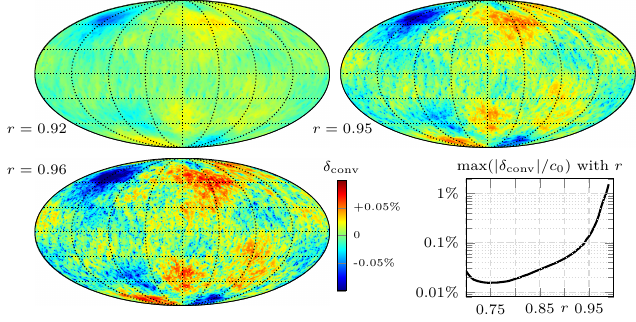}}
  \caption{Wave speed perturbation $\delta_{\mathrm{conv}}$ 
           \cref{eq:wavespeed-convection} generated from solar 
           convection, \cite{Noraz2025}, given as a percentage
           perturbation of the background model $c_0$. 
           The perturbation is non-zero
           for $r\in(0.7,0.99)$, and we show the Mollweide
           map at height \num{0.92} (\textbf{top left}), 
           \num{0.95} (\textbf{top right}) and \num{0.96} (\textbf{bottom left}).
           The relative difference with the background wave speed 
           is maximal near surface (\textbf{bottom right}), with 
           at most about \num{1}\si{\percent} modification of the 
           background $c_0$.
          }
  \label{fig:convection-models-deltac}
\end{figure}

The perturbations induced by active regions are
localized, affecting only a small portion of the solar 
radius and with clear structures (\cref{fig:active-region:model}).
We now consider a different model for the 3D perturbation of wave speed,
using a sound-speed snapshot of the convection simulation from \cite{Noraz2025},
which results in a widespread perturbation. 
The 3D wave speed $c_{\mathrm{conv}}$ is defined from 
the perturbation $\delta_{\mathrm{conv}}$ such that, 
\begin{equation}\label{eq:wavespeed-convection}
  c_{\mathrm{conv}}(r,\theta,\phi) \,=\, 
  c_0(r) \,\left( \, 1 \, +\, \delta_{\mathrm{conv}}(r,\theta,\phi) \, \right) \,.
\end{equation}
In this model, the perturbation of wave speed is between 
scaled radius $r=0.70$ and $r=0.99$. Contrary to the case of active
region, there are no clear structures in the perturbation, 
cf. \cref{fig:convection-models-deltac}.
The maximum amplitude of the perturbation is reached at $r = \num{0.99}$ 
and represents at most only about \num{1}\si{\percent} of the background 
wave speed $c_0$.

\begin{figure}[ht!]\centering
  \subfloat[Source at $(r_s, \theta_s, \phi_s) = (0.99, \num{60}\si{\degree}, \num{135}\si{\degree})$.]{
  \includegraphics[scale=1.0]{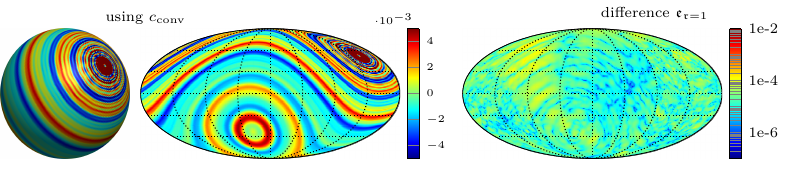}}
  \vspace*{-0.50em}

  \subfloat[Source at $(r_s, \theta_s, \phi_s) = (0.99, \num{-22.5}\si{\degree}, \num{-45}\si{\degree})$.]{
  \includegraphics[scale=1.0]{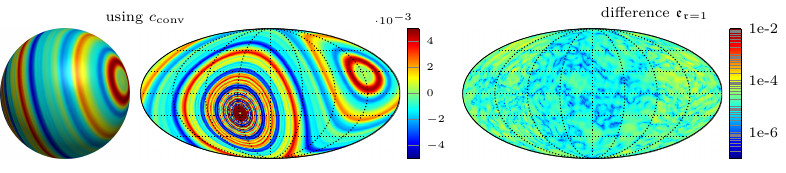}
  }
\caption{Solutions at 1\si{\milli\Hz} using $c_{\mathrm{conv}}$ 
         and relative difference at fixed height $r=1$.}
  \label{fig:convection-wavefields-1mHz}
\end{figure}

In \cref{fig:convection-wavefields-1mHz}, we show the 3D wavefield 
solutions at 1\si{\milli\Hz} using the wave speed $c_{\mathrm{conv}}$ 
for different source positions, alongside with the corresponding differences 
from the solution computed with the radially symmetric background $c_0$. 
The differences vary significantly between the two source locations, 
and no clear pattern corresponding to the perturbation $\delta_{\mathrm{conv}}$ 
can be identified. 
The perturbation $\delta_{\mathrm{conv}}$ affects a wide 
radial interval (for $r \in (0.7, 0.99)$) in which all angles 
are concerned. 
Consequently, no clear pattern in the difference can be observed.

\bigskip

The two background perturbations we have examined highlight 
the influence of both the amplitude and spatial structure of 
the wave speed perturbation on wave propagation. In the case 
of active regions extracted from magnetograms, we have designed 
perturbations that are 
localized and have well-defined spatial patterns. Even with 
large amplitudes, their impact remains confined to these 
locations at 1\si{\milli\Hz} and \num{1.5}\si{\milli\Hz}. 
Increasing the frequency to 2\si{\milli\Hz}, the pattern 
are not clearly visible anymore in the differences.
On the other hand, the perturbations derived from the convection model 
are distributed over a broad radial range, for $r \in (0.7, 0.99)$,
covering all angles. 
Although the amplitude is comparatively small, their coverage causes 
global modifications in the wavefield. 


\section{Conclusion} 
\label{section:conclusion}

In this work, we developed a 3D solver for the time-harmonic 
stellar oscillation equations using the realistic solar background. 
The frequency-domain formulation 
allows us to handle negative buoyancy frequencies, thereby enabling the 
consideration of realistic solar background models including 
three-dimensional heterogeneities.
The discretization is performed using the HDG method which, by defining 
the global matrix in terms of numerical traces only, reduces the size of 
the linear system to be solved compared to traditional DG methods. 
In addition, $hp$-adaptivity is straightforward within the family of 
DG methods, i.e., the polynomial order is adapted on each cell 
to the local wavelength to efficiently control the number of unknowns.

\paragraph{Numerical strategy}

In \cref{section:implementation,section:mumps}, we have established 
the appropriate strategy to carry out the numerical experiments of 3D
gravito-acoustic waves in the Sun.

\begin{enumerate}[leftmargin=*]\setlength{\itemsep}{-1pt}
\item Three formulations of the 
      equations have been introduced, and we have shown that 
      the Liouville variants provide a substantial decrease in the 
      condition number of the matrix. 
      We have used $\Lliouville$ which allows to easily encode 
      3D perturbations of wave speed.
\item The mesh is refined near the surface to account for the decrease 
      in wavelength, and $p$-adaptivity provides flexibility by 
      locally adjusting the polynomial order. 
      Meshes with 6 and 8 millions tetrahedra are used, 
      resulting in matrices of size a few hundred million depending on the frequency.
\end{enumerate}

The computational cost in solving the resulting 
linear system is controlled with \mumps. 
The block low-rank (BLR) compression and the mixed-precision representation of 
the low-rank blocks are employed to reduce both memory usage and runtime, 
and we have investigated the trade-off with the solution accuracy.

\begin{enumerate}[leftmargin=*]\setcounter{enumi}{2}\setlength{\itemsep}{-1pt}
\item The HDG discretization makes block analysis straightforward, thus  
      reducing the computational cost of this phase. In our case, the 
      partitioner \texttt{Scotch} gives the best results.
\item BLR thresholds \blr{7} and \blr{9} 
      are appropriate choices that yield satisfactory 
      accuracy while significantly reducing the factorization 
      cost compared to full-rank. 
      On the other hand, \blr{5} is too aggressive and 
      deteriorates the solution accuracy.
\end{enumerate}
The BLR approximations lead to a decrease in 
both the number of operations and entries in the factors required for 
the factorization compared to the full-rank. 
Only about $10\si{\percent}$ of the operations need to be performed, 
and at most $30\si{\percent}$ of the entries remain.
Furthermore, mixed-precision storage consistently 
contributes to a further $10\si{\percent}$ reduction in the 
storage requirement of the entries.

\paragraph{3D solar wave experiments}

The solver has been validated by comparing with solutions obtained 
under axisymmetric assumptions, highlighting the difficulty of 
accurately capturing the singularity associated with 
near-surface Dirac sources in realistic solar backgrounds. 
This aspect warrants further investigation, for example through 
local mesh refinement or singularity extraction techniques. 
Using Gaussian source functions, we have compared the 
wavefields resulting from spherically symmetric background 
wave speeds with those including 3D perturbations.
We have created localized perturbations representative of 
active regions from observed magnetograms, and have used 
widespread perturbations from solar convection. 

Our results show a linear interdependence between 
perturbations of wave speed and wavefields for localized 
perturbation (such as active region) and relatively low
frequencies. However, increasing the frequency or considering 
widespread perturbation (such as convection) result in the 
absence of clear pattern in the maps of the differences between
wavefield using a spherical background.
This works opens up the perspective of accurate forward modeling that 
can fully adapt to high-resolution observations, taking into account 
three-dimensional heterogeneities and non-spherical stellar structures
and geometries. It enables to simulate the effect of activity on stellar observables and can be used to improve our understanding of farside imaging \cite{Yang2023}. The effect of convection on the
oscillations is also a main topic of research in the stellar community in order to understand 
the differences between the observed and simulated eigenfrequencies (and power spectrum), 
the so-called `surface-effect' \cite{Rosenthal1999}.


\section*{Acknowledgments}

  The authors would like to thank Chiara Puglisi for 
  helpful discussions, and Frank Heckes for 
  his support in using the \mps~cluster.
  This project was provided with computer and storage resources 
  by GENCI at CINES thanks to the grant \texttt{c1715711} on 
  the supercomputer Adastra's GENOA partition.
  This work was partially supported by the ANR-DFG 
  project \textsc{Butterfly}, grant number 
  ANR-23-CE46-0009-01. 
  This work was partially supported by the EXAMA 
  (Methods and Algorithms at Exascale) 
  project under grant ANR--22--EXNU--0002.
  FF acknowledges funding by the European Union 
  with ERC Project \textsc{Incorwave} -- grant 101116288.
  DF and LG acknowledge support from ERC Synergy grant WHOLE SUN 810218.
  Views and opinions expressed are however those of the authors 
  only and do not necessarily reflect those 
  of the European Union or the European Research Council 
  Executive Agency (ERCEA). Neither the European Union nor the
  granting authority can be held responsible for them.
  
\bibliographystyle{siamplain}
\small
\bibliography{sections/bibliography}
\end{document}